\documentclass[11pt,a4]{article}
\usepackage{latexsym}
\usepackage{amssymb,amsmath}
\usepackage{graphics}
\input epsf
\makeatletter
\def\section{\@startsection {section}{1}{\z@}{-3.5ex plus -1ex minus
 -.2ex}{2.3ex plus .2ex}{\large\sc}}
\def\subsection{\@startsection{subsection}{2}{\z@}{-3.25ex plus -1ex minus
 -.2ex}{1.5ex plus .2ex}{\normalsize\sc}}
\makeatother
\makeatletter \@addtoreset{equation}{section}
\renewcommand{\theequation}{\thesection.\arabic{equation}} 
\makeatother
\newcommand{\nc}{\newcommand}

\nc{\bea}{\begin{eqnarray}} \nc{\eea}{\end{eqnarray}}
\nc{\be}{\bea} \nc{\ee}{\eea}
\nc{\tr}{\mathop{\mbox{tr}}\nolimits}
\nc{\ad}{\mathop{\mbox{ad}}\nolimits}
\nc{\Tr}{\mathop{\mbox{Tr}}\nolimits}
\nc{\Det}{\mathop{\mbox{Det}}\nolimits}
\nc{\rk}{\mathop{\mbox{rk}}\nolimits} \nc{\ra}{\rightarrow}
\nc{\Ra}{\Rightarrow} \nc{\LRa}{\Leftrightarrow} \nc{\ot}{\otimes}
\nc{\non}{\nonumber\\} \nc{\ZZ}{\mathbb{Z}} \nc{\RR}{\mathbb{R}}

\newtheorem{theorem}{Theorem}

\newtheorem{lemma}{Lemma}
\newtheorem{corollary}{Corollary}

\def\1#1{^{(#1)}}

\def\ra{\rangle}
\begin{document}
\title{New Identities for Degrees of Syzygies\\ in Numerical Semigroups}
\author{Leonid G. Fel\\
\\Department of Civil Engineering, Technion, Haifa 3200, Israel}
\vspace{-.2cm}
\date{}

\maketitle
\def\be{\begin{equation}}
\def\ee{\end{equation}}
\def\bea{\begin{eqnarray}}
\def\eea{\end{eqnarray}}
\def\p{\prime}
\vspace{-.2cm}
\begin{abstract}
We derive a set of polynomial and quasipolynomial identities for degrees of 
syzygies in the Hilbert series $H\left({\bf d}^m;z\right)$ of nonsymmetric 
numerical semigroups ${\sf S}\left({\bf d}^m\right)$ of arbitrary generating set
of positive integers ${\bf d}^m=\left\{d_1,\ldots,d_m\right\}$, $m\geq 3$. These
identities were obtained by studying together the rational representation of the
Hilbert series $H\left({\bf d}^m;z\right)$ and the quasipolynomial 
representation of the Sylvester waves in the restricted partition function 
$W\!\left(s,{\bf d}^m\right)$. In the cases of symmetric semigroups and complete
intersections these identities become more compact.
\end{abstract}
\section{Introduction}\label{s1}
The study of Diophantine equations is on the border-line between combinatorial 
number theory and commutative algebra. Most important results are awaiting at 
intersection of both theories that has been already seen in the last decades 
\cite{stur05}. 

Focusing on such intersection in study of linear Diophantine equations, in this 
paper we bring together two different approaches, theory of restricted partition
and theory of commutative semigroup rings, and show that such merging is 
fruitful to produce new results. It allows to establish a set of new 
quasipolynomial identities for degrees of the syzygies in the Hilbert series 
$H\left({\bf d}^m;z\right)$ of nonsymmetric numerical semigroups ${\sf S}\left(
{\bf d}^m\right)$ of any embedding dimension, $m\geq 3$, and arbitrary 
generating set of positive integers ${\bf d}^m=\left\{d_1,\ldots,d_m\right\}$ 
where $\gcd(d_1,\ldots,d_m)=1$. The special cases of symmetric semigroups and 
complete intersections make these identities more compact.

Regarding the novelty of these identities, to the best of our knowledge they 
have not been discussed earlier in literature. On the other hand, all necessary 
technical tools to derive them were already elaborated in seminal works of 
Sylvester on partitions (1857, 1882, 1897) and in the basis theorem (1888) and 
the syzygy theorem (1890) of Hilbert.

The paper is organized in six sections. In section \ref{s2} we recall the main 
facts about numerical semigroups ${\sf S}\left({\bf d}^m\right)$ and their 
Hilbert series $H\left({\bf d}^m;z\right)$. For degrees of the syzygies we state
the main result on polynomial and trigonometric identities (Theorem \ref{the1} 
and Corollary \ref{cor1}) which are independent of structure of the generating 
set ${\bf d}^m$. Another result on quasipolynomial identities (Theorem 
\ref{the2}) is valid when among the generators $d_i$ of the set ${\bf d}^m$ 
there exists a subset $\Xi_q\left({\bf d}^m\right)\subset {\bf d}^m$ such that 
$\Xi_q\left({\bf d}^m\right):=\left\{d_i\quad |\quad q\mid d_i\right\}$ and 
$\#\Xi_q\left({\bf d}^m\right)\geq 2$.

In section \ref{s3} we recall the main facts about partition of nonnegative 
integer $s$ into positive integers $\{d_1,\ldots,d_m\}$, each not greater than 
$s$, and about the number of nonequivalent partitions, or representations 
(Reps) $W\!\left(s,{\bf d}^m\right)$. The main emphasis is done on the Sylvester
waves and their symbolic Reps which are useful to find new identities when 
applying the parity claims to the Sylvester waves. 

In section \ref{s4} we derive the quasipolynomial Reps of Sylvester waves with 
trigonometric functions as coefficients and find when their leading terms are 
vanishing (Lemma \ref{lem1}).

In section \ref{s5} we derive the relationship between the rational Rep of the 
Hilbert series $H\left({\bf d}^m;z\right)$ and the quasipolynomial Rep of the 
Sylvester waves in the restricted partition function $W\!\left(s,{\bf d}^m
\right)$. This allows to prove Theorem \ref{the1} (in section \ref{s51}) and 
Theorem \ref{the2} (in section \ref{s53}) on polynomial and quasipolynomial 
identities for degrees of syzygies. 

In section \ref{s6} we discuss different applications of Theorems \ref{the1} and
\ref{the2} to the various kinds of numerical semigroups: complete intersections,
symmetric semigroups (not complete intersections) generated by 4 and 5 elements,
nonsymmetric and pseudosymmetric semigroups generated by 3 elements, and 
semigroups of maximal embedding dimension. We illustrate a validity of 
identities by examples for numerical semigroups which were discussed earlier in 
literature.
\section{Numerical Semigroups and Hilbert Series $H\left({\bf d}^m;z\right)$}
\label{s2}
Throughout the article we assume that the numerical semigroup ${\sf S}\left({\bf
d}^m\right)$ is generated by a minimal set of positive integers ${\bf d}^m=
\left\{d_1,\ldots,d_m\right\}$ with finite complement in ${\mathbb N}$, $\#
\left\{{\mathbb N}\setminus {\sf S}\left({\bf d}^m\right)\right\}<\infty$. We 
study the generating function $H\left({\bf d}^m;z\right)$ of such semigroup 
${\sf S}\left({\bf d}^m\right)$,
\bea
H\left({\bf d}^m;z\right)=\sum_{s\;\in\;{\sf S}\left({\bf d}^m\right)}z^s\;,   
\label{b1}
\eea
which is referred to as {\em the Hilbert series} of ${\sf S}\left({\bf d}^m
\right)$.

Recall the main definitions and facts on numerical semigroups which are
necessary here. A semigroup ${\sf S}\left({\bf d}^m\right)=\left\{s\in{\mathbb N
}\cup\{0\}\;|\;s=\sum_{i=1}^m x_i d_i,\;x_i\in {\mathbb N}\cup\{0\}\right\}$, is
said to be generated by {\em minimal set} of $m$ natural numbers $d_1<\ldots<d_
m$, $\gcd(d_1,\ldots,d_m)=1$, if neither of its elements is linearly 
representable by the rest of elements. It is classically known that $d_1\geq m$
\cite{heku71} where $d_1$ and $m$ are called {\em the multiplicity} and {\em
the embedding dimension (edim)} of the semigroup, respectively. If equality
$d_1=m$ holds then the semigroup ${\sf S}\left({\bf d}^m\right)$ is called of
{\em maximal edim}. {\em The conductor} $c\left({\bf d}^m\right)$ of semigroup 
${\sf S}\left({\bf d}^m\right)$ is defined by $c\left({\bf d}^m\right):=\min
\left\{s\in{\sf S}\left({\bf d}^m\right)\;|\;s+{\mathbb N}\cup\{0\}\subset {\sf 
S}\left({\bf d}^m\right)\right\}$ and related to {\em the Frobenius number} of 
semigroup, $F\left({\bf d}^m\right)=c\left({\bf d}^m\right)-1$. 

A semigroup ${\sf S}\left({\bf d}^m\right)$ is called {\em symmetric} if for 
any integer $s$ the following condition holds: if $s\in {\sf S}\left({\bf d}^m  
\right)$ then $F\left({\bf d}^m\right)-s\not\in {\sf S}\left({\bf d}^m\right)$.
Otherwise ${\sf S}\left({\bf d}^m\right)$ is called nonsymmetric. Notably that
all semigroups ${\sf S}\left(d_1,d_2\right)$ are symmetric.

Denote by $\Delta\left({\bf d}^m\right)$ the complement of ${\sf S}\left({\bf  
d}^m\right)$ in ${\mathbb N}$, i.e., $\Delta\left({\bf d}^m\right)={\mathbb N}
\setminus{\sf S}\left({\bf d}^m\right)$, and call it the set of gaps. The 
cardinality ($\#$) of $\Delta\left({\bf d}^m\right)$ is called {\em the genus}
of ${\sf S}\left({\bf d}^m\right)$, $G\left({\bf d}^m\right):=\#\Delta\left(  
{\bf d}^m\right)$. For the set $\Delta\left({\bf d}^m\right)$ introduce the
generating function $\Phi\left({\bf d}^m;z\right)$ which is related to the
Hilbert series,
\bea
\Phi\left({\bf d}^m;z\right)=\sum_{s\;\in\;\Delta\left({\bf d}^m\right)}z^s\;,
\;\;\;\;\;\;\;\;\Phi\left({\bf d}^m;z\right)+H\left({\bf d}^m;z\right)=
\frac1{1-z}\;.\nonumber
\eea

The Hilbert series $H\left({\bf d}^m;z\right)$ of numerical semigroup ${\sf S}
\left({\bf d}^m\right)$ is a rational function \cite{stan79}
\bea
H\left({\bf d}^m;z\right)=\frac{Q\left({\bf d}^m;z\right)}{\prod_{i=1}^m
\left(1-z^{d_{i}}\right)}\;,\label{tr16}
\eea
where $H\left({\bf d}^m;z\right)$ has a pole $z=1$ of order 1. The numerator  
$Q\left({\bf d}^m;z\right)$ is a polynomial in $z$,
\bea
&&Q\left({\bf d}^m;z\right)=1-Q_1\left({\bf d}^m;z\right)+Q_2\left({\bf d}^m;z
\right)-\ldots+(-1)^{m-1}Q_{m-1}\left({\bf d}^m;z\right)\;,\label{bet05}\\
&&Q_i\left({\bf d}^m;z\right)=\sum_{j=1}^{\beta_i\left({\bf d}^m\right)}z^{C_{
j,i}}\;,\;\;\;1\leq i\leq m-1\;,\;\;\;\deg Q_i\left({\bf d}^m;z\right)<\deg Q_{
i+1}\left({\bf d}^m;z\right)\;.\label{bet1}
\eea
In formula (\ref{bet1}) the positive numbers $C_{j,i}$ and $\beta_i\left({\bf 
d}^m\right)$ denote {\em the degree of the syzygy } and {\em the Betti number}, 
respectively. The latter satisfy the equality \cite{stan79}
\bea
\beta_0\left({\bf d}^m\right)-\beta_1\left({\bf d}^m\right)+\beta_2\left({\bf 
d}^m\right)-\ldots+(-1)^{m-1}\beta_{m-1}\left({\bf d}^m\right)=0\;,\quad\beta_0
\left({\bf d}^m\right)=1\;.\label{bet2}
\eea
The summands $z^{C_{j,i}}$ in (\ref{bet1}) stand for the syzygies of different
kinds and $C_{j,i}$ are the degrees of homogeneous basic invariants for the
syzygies of the $i$th kind,
\bea
&&C_{j,i}\in {\mathbb N}\;,\;\;\;\;C_{j+1,i}\geq C_{j,i}\;,\;\;\;\;C_{\beta_{i+
1},i+1}>C_{\beta_i,i}\;,\;\;\;\;C_{1,i+1}>C_{1,i}\;,\;\;\;\;\mbox{and}
\nonumber\\
&&C_{j,i}\neq C_{r,i+2k-1}\;,\;\;\;1\leq j\leq\beta_i\left({\bf d}^m\right)\;,
\;\;\;1\leq r\leq \beta_{i+2k-1}\left({\bf d}^m\right)\;,\;\;\;
1\leq k\leq \left\lfloor \frac{m-i}{2}\right\rfloor\;.\;\;\;\label{ar10a}
\eea
The last requirement (\ref{ar10a}) means that all necessary cancellations  
(annihilations) of terms $z^{C_{j,i}}$ in (\ref{bet05}) are already performed.
However the other equalities, $C_{j,i}=C_{r,i+2k}$ and $C_{j,i}=C_{q,i}$, $j\neq
q$, are not forbidden excluding the syzygy degrees of the last $(m-1)$th kind  
(\cite{fa09}, Lemma 2). 

A degree of the polynomial $Q\left({\bf d}^m;z\right)$ is strongly related 
\cite{fa09} to the Frobenius number by
\bea
\deg Q\left({\bf d}^m;z\right)=F\left({\bf d}^m\right)+\sigma_1\;,\;\;\;\;
\mbox{where}\;\;\;\;\sigma_1=d_1+\ldots +d_m\;.\label{bet1a}
\eea
We present here the Hilbert series for nonsymmetric semigroup generated by 
triple $\{3,5,7\}$,
\bea
H\left(\{3,5,7\};z\right)=\frac{1-z^{10}-z^{12}-z^{14}+z^{17}+z^{19}}{
\left(1-z^3\right)\left(1-z^5\right)\left(1-z^7\right)}\;.\label{ar10b}
\eea

For the rigorous notions of syzygies of the 1st and higher kinds, their moduli
and specifying homomorphisms as well as the Betti numbers and minimal free 
resolution we refer to the book \cite{eis05}. An informal description of 
syzygies, difference binomials and other homogeneous bases for the higher 
syzygies, which came by applying the Hilbert basis and syzygy theorems, can be 
found in the review \cite{stan79}. Regarding the degrees of the syzygies, in the
general case of nonsymmetric numerical semigroups ${\sf S}\left({\bf d}^m
\right)$, their values $C_{j,i}$ are usually obtained by computational algorithm
\cite{del08} calculating the homogeneous bases for the syzygies moduli and their
specifying homomorphisms in a minimal free resolution. 

If a semigroup ${\sf S}\left({\bf d}^m\right)$ is symmetric then a duality 
relation for numerator $Q\left({\bf d}^m;z\right)$ holds \cite{fa09}
\bea
Q\left({\bf d}^m;\frac{1}{z}\right)z^{\deg Q({\bf d}^m;z)}=(-1)^{m-1}
Q\left({\bf d}^m;z\right)\;,\label{teb2}
\eea
and by consequence of (\ref{teb2}) we have
\bea  
\beta_k\left({\bf d}^m\right)=\beta_{m-k-1}\left({\bf d}^m\right)\;,\quad 
C_{j,k}+C_{j,m-k-1}=\deg Q\left({\bf d}^m;z\right)\;,\quad \left\{\begin{array}
{c}0\leq k<m,\\1\leq j\leq\beta_k\left({\bf d}^m\right)\;.\end{array}\right.
\label{teb3}
\eea
In fact, the 2nd equality in (\ref{teb3}) does not contribute much for 
determination of the degrees of the syzygies, since in accordance with 
(\ref{bet1a}) the degree of $Q\left({\bf d}^m;z\right)$ is strongly related to 
the Frobenius number $F\left({\bf d}^m\right)$ which is unknown for $m\geq 3$ 
in terms of generators $d_i$ only.

If $\beta_1\left({\bf d}^m\right)=m-1$ then the semigroup ${\sf S}\left({\bf d}
^m\right)$ is called {\em complete intersection} and the numerator of the 
corresponding Hilbert series is given by \cite{stan79}
\bea
Q\left({\bf d}^m;z\right)=\prod_{j=1}^{m-1}\left(1-z^{e_j}\right)\;,\;\;\;
e_j\in{\mathbb N}\;.\label{teb4}
\eea

The degrees of syzygies $C_{j,i}$ provide more precise and accurate 
characteristics of numerical semigroups than the Betti numbers $\beta_i\left(
{\bf d}^m\right)$. There are only few sorts of semigroups where both sets of 
$\beta_i\left({\bf d}^m\right)$ and $C_{j,i}$ are known completely. This is 
true, e.g., for semigroups of maximal edim (\cite{sal79}, Theorem 1 and 
\cite{f10}, sect. 7) and the ${\sf S}\left({\bf d}^3\right)$ semigroups of 
special kind: Pythagorean semigroups (\cite{fel04}, sect. 6.1), pseudosymmetric 
semigroups (\cite{f10}, sect. 6.1) and Fibonacci and Lucas symmetric semigroups 
\cite{fel091}. A family of semigroups with known Betti numbers but unknown 
degrees of syzygies is much wider, e.g., nonsymmetric semigroups ${\sf S}\left(
{\bf d}^3\right)$ \cite{herz70}, symmetric semigroups ${\sf S}\left({\bf d}^4
\right)$ \cite{ber752}, symmetric semigroups of almost maximal edim, $d_1=m+1$ 
(\cite{sal79}, Theorem 2) and complete intersections \cite{stan79}.
\subsection{Main results}\label{s21}
In this section we present two theorems on polynomial and quasipolynomial 
identities for the degrees $C_{j,i}$ of the syzygies. Their proof will follow 
later, Theorem \ref{the1} in section \ref{s51} and Theorem \ref{the2} in section
\ref{s53}. Start with a generic case of the generating set ${\bf d}^m$ not 
keeping in mind any relationships among the generators $d_j$.
\begin{theorem}\label{the1}
Let the numerical semigroup ${\sf S}\left({\bf d}^m\right)$ be given with its
Hilbert series $H\left({\bf d}^m;z\right)$ in accordance with (\ref{tr16}) and 
(\ref{bet05}). Then the following polynomial identities hold,
\bea
&&\sum_{j=1}^{\beta_1\left({\bf d}^m\right)}C_{j,1}^r-\sum_{j=1}^{\beta_2\left(
{\bf d}^m\right)}C_{j,2}^r+\ldots+(-1)^m\sum_{j=1}^{\beta_{m-1}\left({\bf d}^m
\right)}C_{j,m-1}^r=0\;,\;\;\;\;\;r=1,\ldots ,m-2\;,\quad\label{10c}\\
&&\sum_{j=1}^{\beta_1\left({\bf d}^m\right)}C_{j,1}^{m-1}-\sum_{j=1}^{\beta_2
\left({\bf d}^m\right)}C_{j,2}^{m-1}+\ldots+(-1)^m\sum_{j=1}^{\beta_{m-1}\left(
{\bf d}^m\right)}C_{j,m-1}^{m-1}=(-1)^m(m-1)!\prod_{i=1}^md_i\;.\hspace{1.5cm}
\label{10d}
\eea
\end{theorem}
Additional sort of identities appear when we keep in mind relationships among
the generators $d_i\in {\bf d}^m$. Namely, if there exists a subset $\Xi_q\left(
{\bf d}^m\right)\subset {\bf d}^m$ such that $\Xi_q\left({\bf d}^m\right):=
\left\{d_i\quad |\quad q\mid d_i\right\}$, $\omega_q=\#\Xi_q\left({\bf d}^m
\right)$, then there hold another type of identities which are not polynomials.
\begin{theorem}\label{the2}
Let the numerical semigroup ${\sf S}\left({\bf d}^m\right)$ be given with its
Hilbert series $H\left({\bf d}^m;z\right)$ in accordance with (\ref{tr16}) and
(\ref{bet05}). Then for every $1<q\leq\max\left\{d_1,\ldots,d_m\right\}$, and 
$\gcd(n,q)=1$, $1\leq n<q/2$, the following quasipolynomial identities hold,
{\footnotesize
\bea
\sum_{j=1}^{\beta_1\left({\bf d}^m\right)}\!C_{j,1}^r\exp\left(i\frac{2\pi n}{
q}C_{j,1}\right)-\!\sum_{j=1}^{\beta_2\left({\bf d}^m\right)}\!C_{j,2}^r\exp
\left(i\frac{2\pi n}{q}C_{j,2}\right)+\ldots+(-1)^m\!\sum_{j=1}^{\beta_{m-1}
\left({\bf d}^m\right)}\!C_{j,m-1}^r\exp\left(i\frac{2\pi n}{q}C_{j,m-1}\right)
\!=0\nonumber
\eea}
where $r=1,\ldots ,\omega_q-1$. However, in the case $r=0$ another trigonometric
identity holds,
{\footnotesize
\bea
\sum_{j=1}^{\beta_1\left({\bf d}^m\right)}\!\exp\left(i\frac{2\pi n}{q}C_{j,1}
\right)-\!\sum_{j=1}^{\beta_2\left({\bf d}^m\right)}\!\exp\left(i\frac{2\pi n}{
q}C_{j,2}\right)+\ldots+(-1)^m\!\sum_{j=1}^{\beta_{m-1}\left({\bf d}^m
\right)}\!\exp\left(i\frac{2\pi n}{q}C_{j,m-1}\right)\!=1\;.\label{10i}
\eea}
\end{theorem}
By Theorem \ref{the2} another statement comes irrespectively to the inner
relationships between the generators $d_i$. Indeed, by consequence of 
(\ref{10i}) and the fact that the generating set ${\bf d}^m$ is minimal we have 
Corollary.
\begin{corollary}\label{cor1}
Let the numerical semigroup ${\sf S}\left({\bf d}^m\right)$ be given with its
Hilbert series $H\left({\bf d}^m;z\right)$ in accordance with (\ref{tr16}) and  
(\ref{bet05}). Then for every $d_k\in {\bf d}^m$, $1\leq k\leq m$, and $\gcd(
n,d_k)=1$, $1\leq n<d_k$, the following trigonometric identities hold,
{\footnotesize
\bea
\sum_{j=1}^{\beta_1\left({\bf d}^m\right)}\!\exp\left(i\frac{2\pi n}{d_k}C_{j,1}
\right)-\!\sum_{j=1}^{\beta_2\left({\bf d}^m\right)}\!\exp\left(i\frac{2\pi n}{
d_k}C_{j,2}\right)+\ldots+(-1)^m\!\sum_{j=1}^{\beta_{m-1}\left({\bf d}^m
\right)}\!\exp\left(i\frac{2\pi n}{d_k}C_{j,m-1}\right)\!=1\;.\label{10e}
\eea}
\end{corollary}
In section \ref{s6} we will discuss more special cases of numerical semigroups 
when a part of identities (\ref{10c}) and (\ref{10d}) becomes trivial (see 
Corollary \ref{cor2} for complete intersections) or even they all do not provide
new relations (see section \ref{s611} for telescopic semigroups).
\section{Restricted Partition Functions $W\!\left(s,{\bf d}^m\right)$}\label{s3}
The restricted partition function $W\!\left(s,{\bf d}^m\right)$ is a number of 
partitions of $s$ into positive integers  $\{d_1,\ldots,d_m\}$, each not greater
than $s$. The generating function $M\!\left({\bf d}^m;z\right)$ for $W\!\left(
s,{\bf d}^m\right)$ has a form \cite{and76},
\be
M\!\left({\bf d}^m;z\right)=\prod_{i=1}^m\frac{1}{1-z^{d_{i}}}=\sum_{s=0}^{
\infty}W\!\left(s,{\bf d}^m\right)\;z^s\;,\label{g1}
\ee
where $W\!\left(s,{\bf d}^m\right)$ satisfies the basic recursive relation
\be
W\!\left(s,{\bf d}^m\right)-W\!\left(s-d_k,{\bf d}^m\right)=W\!\left(s,{\bf d}_
k^{m-1}\right),\;\;\;\;\;\;{\bf d}_k^{m-1}=\{d_1,\ldots,d_{k-1},d_{k+1},\ldots,
d_m\}\;.\label{g2}
\ee
The function $W\!\left(s,{\bf d}^m\right)$ is also satisfied parity properties 
\cite{fr02}, Lemma 4.1,
\bea
W\!\left(s-\frac{\sigma_1}{2},{\bf d}^{2m}\right)=-W\!\left(-s-\frac{\sigma_1}{
2},{\bf d}^{2m}\right),\quad W\!\left(s-\frac{\sigma_1}{2},{\bf d}^{2m+1}\right)
=W\!\left(-s-\frac{\sigma_1}{2},{\bf d}^{2m+1}\right),\quad\label{g2a}
\eea
and the statement (\cite{fr02}, Lemma 4.3) about (not all) zeroes $\varsigma_0
\left({\bf d}^m\right)$ of $W\!\left(s,{\bf d}^m\right)$: if the generators 
$d_i$ are mutually prime numbers, $\gcd(d_i,d_k)=\delta_{ik}$, then 
\bea
\varsigma_0\left({\bf d}^{2m+1}\right)=-1,-2,\ldots,-\sigma_1+1\;\;\;\;\mbox{
and}\;\;\;\;\varsigma_0\left({\bf d}^{2m}\right)=-1,-2,\ldots,-\sigma_1+1,
-\frac{\sigma_1}{2}\;.\label{g2b}
\eea
According to definitions of $W\!\left(s,{\bf d}^m\right)$ and $F\left({\bf d}^m
\right)$ of semigroup ${\sf S}\left({\bf d}^m\right)$, given in section \ref{s2},
the Frobenius number is a maximal zero of $W\!\left(s,{\bf d}^m\right)$.

Regarding the generating function $M\left({\bf d}^m;z\right)$ in (\ref{g1}), the
degree of the numerator in $z$ vanishes ($\deg_z 1=0$) while the degree of the 
denominator in $z$ is positive, $\deg_z\prod_{i=1}^m\left(1-z^{d_{i}}\right)=
\sigma_1$. Moreover, every zero $\zeta$ of the denominator satisfies $\zeta^{T}
=1$, where $T=lcm\left({\bf d}^m\right)$ is a least common multiple of the set 
$\{d_1,\ldots,d_m\}$. Then by \cite{stan79}, Proposition 4.4.1, the function 
$W\!\left(s,{\bf d}^m\right)$ is a quasipolynomial of degree $m-1$, 
\bea
W\!\left(s,{\bf d}^m\right)=K_1\left(s,{\bf d}^m\right)\;s^{m-1}+K_2\left(s,
{\bf d}^m\right)s^{m-2}+\ldots +K_{m-1}\left(s,{\bf d}^m\right)s+K_m\left(s,
{\bf d}^m\right)\;,\label{g3}
\eea
where each $K_j\left(s,{\bf d}^m\right)$ is a periodic function with integer 
period $\tau_j$ dividing $T$. According to  Schur's theorem (see \cite{wil06}, 
Theorem 3.15.2) the 1st coefficient $K_1\left(s,{\bf d}^m\right)$ is independent
of $s$,
\bea
K_1\left(s,{\bf d}^m\right)=\frac{1}{(m-1)!\;\pi_m}\;,\;\;\;\pi_m=\prod_{i=1}^m
d_i\;.\label{g4}
\eea
As for the other $K_j\left(s,{\bf d}^m\right)$, $1<j\leq m$, they can be 
calculated by computational algorithm \cite{fr02}, Appendix A, based on 
recursion relation (\ref{g2}) and zeroes' sequences (\ref{g2b}).

The partition function $W\!\left(s,{\bf d}^m\right)$ for the set of consecutive 
integers $\{d_1=1,\ldots,d_m=m\}=\{\overline{m}\}$ is called {\em unrestricted}.
It was under special consideration in \cite{fr02} and the corresponding formulas
for $W\!\left(s,\{\overline{m}\}\right)$, $m\leq 12$, were also presented there.
Although a variable $s$ in formulas (\ref{g1}) -- (\ref{g3}) is assumed to have 
integer values, but Rep (\ref{g3}) can be extended to real values of $s$, though
such extension is not unique. It does depend on the choice of complete set of 
periodic functions for $K_j\left(s,{\bf d}^m\right)$ which is often taken as 
${\sf sin}\left(\frac{2\pi k}{T}s\right)$ and ${\sf cos}\left(\frac{2\pi k}{T}
s\right)$, $k\in{\mathbb N}$, e.g., \cite{kom03}
\bea
W\!\left(s,\{3,5,7\}\right)=\frac{s^2}{210}+\frac{s}{14}+\frac{74}{315} +\frac{
2}{9}\cos\frac{2\pi s}{3}+\frac{8}{25}\left[\left(\sin\frac{\pi}{5}\right)^2\!\!
\cos\frac{2\pi s}{5}+\left(\sin\frac{2\pi}{5}\right)^2\!\!\cos\frac{4\pi s}{5}
\right]\nonumber\\
-\frac{2}{7\sqrt{7}}\left[\sin\frac{6\pi}{7}\cos\frac{6\pi s}{7}+2\left(\sin 
\frac{\pi}{7}\right)^2\sin\frac{2\pi s}{7}+2\left(\sin\frac{2\pi}{7}\right)^2
\sin\frac{4\pi s}{7}\right]\nonumber\\
+\frac{2}{7\sqrt{7}}\left[\sin\frac{2\pi}{7}\cos\frac{2\pi s}{7}+\sin\frac{4\pi}
{7}\cos\frac{4\pi s}{7}+2\left(\sin\frac{3\pi}{7}\right)^2\sin\frac{6\pi s}{7}
\right].\quad\label{g5a}
\eea
A brief comparison of the last formula with (\ref{g3}) shows that $K_2\left(s,\{
3,5,7\}\right)$ takes a constant value, this is much stronger than to be 
periodic in $s$ with period dividing $lcm(3,5,7)=120$. Similar phenomenon holds 
for other generating sets ${\bf d}^m$, $m\geq 3$, with mutually prime generators
$\gcd(d_i,d_k)=\delta_{ik}$, i.e., $K_j\left(s,{\bf d}^m\right)$, $j=1,\ldots,m-
1$, does not depend on $s$ \cite{rf06}. This indicates that the basic properties
(\ref{g2}), (\ref{g2a}) and (\ref{g2b}) of $W\!\left(s,{\bf d}^m\right)$ and 
general considerations regarding the generating function $M\!\left({\bf d}^m;z
\right)$ in (\ref{g1}), which were preceding a quasipolynomial Rep (\ref{g3}), 
are quite weak to study $K_j\left(s,{\bf d}^m\right)$ in more details.
\subsection{Sylvester Waves}\label{s31}
A powerful approach to study $W\!\left(s,{\bf d}^m\right)$ dates back to 
Sylvester \cite{sy97} and his recipe enabling to determine a restricted 
partition function by decomposing it into {\em Sylvester waves} $W_q\!\left(s,
{\bf d}^m\right)$,
\bea
W\!\left(s,{\bf d}^m\right)=\sum_{q=1,\;q\mid d_i}^{\max{\bf d}^m} W_q\!\left(
s,{\bf d}^m\right)\;,\quad \max{\bf d}^m=\max\{d_1,\ldots,d_m\}\;,\label{g6}
\eea
where summation runs over all distinct factors of $m$ generators $d_i$. By 
(\ref{g6}) every wave $W_q\!\left(s,{\bf d}^m\right)$ is a quasipolynomial in 
$s$ and by \cite{rf06}, section 5, it satisfies the recursive relation 
(\ref{g2}).

Sylvester stated and proved \cite{sy97} that the wave $W_q\!\left(s,{\bf d}^m
\right)$ is a residue at a point $z=0$ of a function $F_q(s,z)$ which is given 
by
\bea
F_q(s,z)=\sum_{1\leq n<q\atop\gcd(n,q)=1}\frac{\xi_q^{-s\;n}\;e^{s\;z}}{\prod_{
k=1}^{m}\left(1-\xi_q^{d_k n}\;e^{-d_k z}\right)}\;,\;\;\quad\;\xi_q=\exp\left(
\frac{2\pi i}{q}\right)\;.\label{g7}
\eea
The summation in (\ref{g7}) is made over all prime roots $\xi_q^n$ for $n$ 
relatively prime to $q$ (including unity) and smaller than $q$. Making use of 
Rep (\ref{g7}) Sylvester showed that every wave $W_q\!\left(s,{\bf d}^m\right)$ 
possesses also the parity property (\ref{g2a}),
\bea
W_q\!\left(s-\frac{\sigma_1}{2},{\bf d}^{2m}\right)=-W_q\!\left(-s-\frac{\sigma
_1}{2},{\bf d}^{2m}\right),\;W_q\!\left(s-\frac{\sigma_1}{2},{\bf d}^{2m+1}
\right)=W_q\!\left(-s-\frac{\sigma_1}{2},{\bf d}^{2m+1}\right)\;\;\label{g10c}
\eea

The waves $W_q\!\left(s,{\bf d}^m\right)$ were found \cite{rf06} in a form of 
finite sum of the Bernoulli polynomials of higher order multiplied by $q$ - 
periodic function expressed through the Eulerian polynomials of higher order. In
this section we give symbolic formulas for $W_q\!\left(s,{\bf d}^m\right)$ which
are more appropriate when dealing with higher $q$ \cite{rf06}. The 1st wave $W_1
\!\left(s,{\bf d}^m\right)$ is a polynomial part of the whole $W\!\left(s,{\bf 
d}^m\right)$ and serves as a good approximant for the whole $W\left(s,{\bf d}^m
\right)$ \cite{rf06},
\bea
W_1\!\left(s,{\bf d}^m\right)=\frac{1}{(m-1)!\;\pi_m}\left(s+\sigma_1+\sum_{i=
1}^m B\;d_i\right)^{m-1}\;.\label{g8}
\eea
As it is convenient in symbolic (umbral) calculus \cite{ro84}, in (\ref{g8}) 
after binomial expansion the powers $(B \;d_i)^r$ are converted into the 
generator's powers multiplied by Bernoulli numbers, i.e., $d_i^rB_r$. More 
details will be given in section \ref{s41}.

The 2nd wave $W_2\!\left(s,{\bf d}^m\right)$ reads in symbolic form \cite{rf06},
\bea
W_2\!\left(s,{\bf d}^m\right)=\frac{2^{\omega_2-m}\cos \pi s}{(\omega_2-1)!\;
\pi_{\omega_2}}\left(s+\sigma_1+\sum_{i=1}^{\omega_2}B\;d_i+\sum_{i=\omega_2+1}
^{m}E(0)\;d_i\right)^{\omega_2-1},\label{g9}
\eea
where $\omega_2$ and $\pi_{\omega_2}$ are related to the set $\Xi_2\left({\bf 
d}^m\right)$ comprising only the even generators $d_i$,
\bea
\Xi_2\left({\bf d}^m\right):=\left\{d_i\quad |\quad 2\mid d_i\right\},\quad
\omega_2=\#\Xi_2\left({\bf d}^m\right),\quad \pi_{\omega_2}=\prod_{d_i\in\Xi_2
\left({\bf d}^m\right)}d_i\;.\nonumber
\eea
As in formula (\ref{g8}), the symbolic binomial expansion the powers $\left(E
(0)\;d_i\right)^r$ in (\ref{g9}) has to be converted into the generator's powers
multiplied by the values of the Euler polynomial $E_r(x)$ at $x=0$, i.e., $d_i^r
E_r(0)$. Note that $E_r(0)$ differs from the Euler number $E_r=2^rE_r(1/2)$.

The $q$-th wave $W_q\!\left(s,{\bf d}^m\right)$, $q>1$, reads in symbolic 
form \cite{rf06},
\bea
W_q\!\left(s,{\bf d}^m\right)=\frac{1}{(\omega_q-1)!\;\pi_{\omega_q}}\sum_{1\leq
n<q\atop\gcd(n,q)=1}\!\!\!\!{\cal W}_{q,n}\left(s,{\bf d}^m\right),\quad\quad
\mbox{where}\label{g10}
\eea
\bea
{\cal W}_{q,n}\!\left(s,{\bf d}^m\right)=\frac{\xi_q^{-s\;n}}{\prod_{i=\omega_q+
1}^m\left(1-\xi_q^{d_i\;n}\right)}\left(s+\sigma_1+\sum_{i=1}^{\omega_q}B\;d_i
+\!\!\!\sum_{i=\omega_q+1}^{m}\!\!\!H\left(\xi_q^{d_i\;n}\right)\;d_i\right)^{
\omega_q-1},\label{g10e}
\eea
and $\omega_q$ and $\pi_{\omega_q}$ are related to the set $\Xi_q\left({\bf d}^m
\right)$ comprising only the generators $d_i$ divided by $q$,
\bea
\Xi_q\left({\bf d}^m\right):=\left\{d_i\quad |\quad q\mid d_i\right\},\quad
\omega_q=\#\Xi_q\left({\bf d}^m\right),\quad\pi_{\omega_q}=\prod_{d_i\in\Xi_q
\left({\bf d}^m\right)}d_i\;.\nonumber
\eea
The numbers $\left(H\left(\xi_q^{d_i\;n}\right)\right)^r=H_r\left(\xi_q^{d_i\;n}
\right)$ generalize the corresponding $E_r(0)=H_r(-1)$. They were introduced by 
Frobenius \cite{Fr10} and Carlitz \cite{Ca60} as the values of the rational 
function $H_n(x)$ at $x=\xi_q^{d_i\;n}$, where $H_n(x)$ itself comes by power 
expansion of its generating function,
\bea
\frac{1-x}{e^t-x}=1+\sum_{n=1}^{\infty}H_n(x)\frac{t^n}{n!},\quad\quad\mbox{
where}\quad H_n\left(x^{-1}\right)\!=\!(-1)^nx\;H_n(x)\;,\quad H_{2n}(-1)=0\;,
\label{g10b}
\eea
and for $x\neq 1$ the rational function $H_n(x)$ read
\bea
H_1(x)=\frac1{x-1},\;\;H_2(x)=\frac{x+1}{(x-1)^2},\;\;H_3(x)=\frac{x^2+4x+1}{
(x-1)^3},\;\;H_4(x)=\frac{x^3+11x^2+11x+1}{(x-1)^4},\;
\ldots\nonumber
\eea

Consider a special case of the tuple ${\bf p}^m=\{p_1,p_2,\ldots,p_m\}$ of 
primes $p_i$ which leads to essential simplification of formula (\ref{g6}). The 
1st Sylvester wave is given by (\ref{g8}) while all higher waves are purely 
periodic \cite{rf06},
\bea
W_{p_{i}}\left(s;{\bf p}^m\right)=\frac{1}{p_i}\sum_{n=1}^{p_{i}-1}\frac{
\xi_{p_{i}}^{-sn}}{\prod_{r\neq i}^m\left(1-\xi_{p_{i}}^{p_r n}\right)}\;.
\label{g11}
\eea
Calculating $W_q\left(s;\{3,5,7\}\right)$, $q=1,3,5,7$, one can get (\ref{g5a})
that explains why $K_2\left(s,\{3,5,7\}\right)$ in formula (\ref{g3}) is taking 
a constant value. We arrive at the similar conclusion for the generating sets 
${\bf d}^m$, $m\geq 3$, $\gcd(d_i,d_k)=\delta_{ik}$ : since $\omega_q=1$ in 
(\ref{g10}) then coefficients $K_j\left(s,{\bf d}^m\right)$, $j=1,\ldots, m-1$, 
in (\ref{g3}) are taking constant values. 
\section{Quasipolynomial Representation for $W_q\!\left(s,{\bf d}^m\right)$}
\label{s4}
In this section we specify the quasipolynomial Reps of the Sylvester waves with 
trigonometric functions as coefficients. 
Start with technical details and note that by identity $\xi_q^{n_1}=\xi_q^{  
-n_2}$, $n_1+n_2=q$, a sum of two partial Sylvester waves ${\cal W}_{q,n}\left(
s,{\bf d}^m\right)$ and ${\cal W}_{q,q-n}\left(s,{\bf d}^m\right)$ in
(\ref{g10}) can be represented as follows
\bea
\overline{{\cal W}}_{q,n}\left(s,{\bf d}^m\right)={\cal W}_{q,n}\left(s,{\bf d}
^m\right)+{\cal W}_{q,-n}\left(s,{\bf d}^m\right)\;.\label{g10o}
\eea
Thus, instead of (\ref{g10}) we write
\bea
W_q\!\left(s,{\bf d}^m\right)=\frac{1}{(\omega_q-1)!\;\pi_{\omega_q}}\sum_{1
\leq n<q/2\atop\gcd(n,q)=1}\!\!\!\!\overline{{\cal W}}_{q,n}\left(s,{\bf d}^m
\right)\;.\label{g10r}
\eea
Explicit formulas (\ref{g10e}), (\ref{g10o}) and (\ref{g10r}) for Sylvester 
waves have one serious lack: their expressions are highly cumbersome and 
difficult to deal with. On the other hand, a visible simplicity of 
quasipolynomial (\ref{g3}) is accompanied by another lack: periodic functions
$K_j\left(s,{\bf d}^m\right)$ don't distinguish between harmonics with distinct
periods. For the purpose of this article it would be worthwhile to have 
something intermediate, rather simple but still inherited basic properties of 
(\ref{g10}) and (\ref{g10e}) even if a minor portion of information would left 
unknown.

Keeping in mind (\ref{g10e}) choose the following Rep for ${\cal W}_{q,\pm n}
\left(s,{\bf d}^m\right)$,
\bea
{\cal W}_{q,\pm n}\left(s,{\bf d}^m\right)=\frac{{\sf L}^{q,n}\left(s,{\bf d}^m
\right)\pm i\;{\sf M}^{q,n}\left(s,{\bf d}^m\right)}{2}\cdot\xi_q^{\mp sn}\;,
\label{g10z}
\eea
where ${\sf L}^{q,n}_m(s)$ and ${\sf M}^{q,n}_m(s)$ are real functions. 
Inserting $\xi_q$ from (\ref{g7}) into (\ref{g10z}) and (\ref{g10o}) we get
\bea
\overline{{\cal W}}_{q,n}\left(s,{\bf d}^m\right)={\sf L}^{q,n}\left(s,{\bf d}^m
\right)\cdot\cos\frac{2\pi n}{q}s+{\sf M}^{q,n}\left(s,{\bf d}^m\right)\cdot\sin
\frac{2\pi n}{q}s\;.\label{g10y}
\eea
Formulas (\ref{g3}) and (\ref{g6}) for restricted partition function allow to 
construct one more Rep that reflects the basic properties of the both of them: 
$W\!\left(s,{\bf d}^m\right)$ is a real function comprising the quasimonomial 
terms $s^k\sin\left(\frac{2\pi n}{q}s\right)$ and $s^k\cos\left(\frac{2\pi n}{
q}s\right)$, $1\leq k\leq m$. Thus, 
\bea
W_1\!\left(s,{\bf d}^m\right)=\frac{1}{(m-1)!\;\pi_m}\;{\sf L}^1\left(s,{\bf d}
^m\right),\quad W_2\!\left(s,{\bf d}^m\right)=\frac{2^{\omega_2-m}}{
(\omega_2-1)!\;\pi_{\omega_2}}\;{\sf L}^{2,1}\left(s,{\bf d}^m\right)\cdot\cos
\pi s\;,\label{g10t}
\eea
while $W_q\!\left(s,{\bf d}^m\right)$ and $\overline{{\cal W}}_{q,n}\!\left(s,
{\bf d}^m\right)$, $q\geq 3$, are given in (\ref{g10r}) and (\ref{g10y}), 
respectively. Polynomials ${\sf L}^1\left(s,{\bf d}^m\right)$, ${\sf L}^{q,n}
\left(s,{\bf d}^m\right)$ and ${\sf M}^{q,n}\left(s,{\bf d}^m\right)$ read
\bea
{\sf L}^1\left(s,{\bf d}^m\right)&=&\ell_1\left({\bf d}^m\right)s^{m-1}+
\ell_2\left({\bf d}^m\right)s^{m-2}+\ldots+\ell_{m-1}\left({\bf d}^m\right)s+
\ell_{m}\left({\bf d}^m\right),\label{g12a}\\
{\sf L}^{q,n}\left(s,{\bf d}^m\right)&=&L_1^{q,n}\left({\bf d}^m\right)
s^{\omega_q-1}+L_2^{q,n}\left({\bf d}^m\right)s^{\omega_q-2}+\ldots+
L_{\omega_q}^{q,n}\left({\bf d}^m\right),\quad q\geq 2,\label{g13}\\
{\sf M}^{q,n}\left(s,{\bf d}^m\right)&=&M_1^{q,n}\left({\bf d}^m\right)
s^{\omega_q-1}+M_2^{q,n}\left({\bf d}^m\right)s^{\omega_q-2}+\ldots+
M_{\omega_q}^{q,n}\left({\bf d}^m\right),\quad q\geq 3.\label{g13a}
\eea
Coefficients $\ell_j\left({\bf d}^m\right)$, $L_j^{q,n}\left({\bf d}^m\right)$ 
and $M_j^{q,n}\left({\bf d}^m\right)$ in (\ref{g12a}) -- (\ref{g13a}) are real 
and so far unknown. 
\subsection{The 1st Sylvester wave}\label{s41}
Performing a binomial expansion in (\ref{g8}) and comparing it with (\ref{g10t})
we get
\bea
{\sf L}^1\left(s,{\bf d}^m\right)=\sum_{r=0}^{m-1}{m-1\choose r}f_r\cdot s^{m-
1-r}\;,\quad\quad f_r=\left(\sigma_1+\sum_{i=1}^mB\;d_i\right)^r\;,\label{e1}
\eea
where for $f_r$ the above formula presumes a symbolic exponentiation. Denoting 
$\sigma_k=\sum_{i=1}^md_i^k$ and making use of the sequence of Bernoulli 
numbers,
\bea
B_0=1,\;B_1=-\frac{1}{2},\;B_2=\frac{1}{6},\;B_4=-\frac{1}{30},\;B_6=\frac{1}
{42}\;,\ldots,\;\;B_{2k+1}=0\;,\;k\geq 1\;,\label{e1a}
\eea
we perform a straightforward calculation in (\ref{e1}) and give the seven first 
values $f_r$, $0\leq r\leq 6$,
\bea 
&&f_0=1\;,\quad f_1=\frac{\sigma_1}{2}\;,\quad f_2=\frac{1}{2^2}\left(\sigma_1^2
-\frac{\sigma_2}{3}\right)\;,\quad f_3=\frac{\sigma_1}{2^3}\left(\sigma_1^2-
\sigma_2\right)\;,\label{e2}\\
&&f_4=\frac{1}{2^4}\left(\sigma_1^4-2\sigma_1^2\sigma_2+\frac{\sigma_2^2}{3}+
\frac{2\sigma_4}{15}\right)\;,\quad f_5=\frac{\sigma_1}{2^5}\left(\sigma_1^4-
\frac{10\sigma_1^2\sigma_2}{3}+\frac{5\sigma_2^2}{3}+\frac{2\sigma_4}{3}
\right)\;,\nonumber\\
&&f_6=\frac1{2^6}\left(\sigma_1^6-5\sigma_1^4\sigma_2+5\sigma_1^2\sigma_2^2+
2\sigma_1^2\sigma_4-\frac{5\sigma_2^3}{9}-\frac{2\sigma_2\sigma_4}{3}-
\frac{16\sigma_6}{63}\right)\;.\nonumber
\eea
In fact, Sylvester's paper \cite{sy97} contains already formulas for $W_1\!
\left(s,{\bf d}^m\right)$, $1\leq m\leq 7$, where one can recognize $f_r$ given 
in (\ref{e2}). By comparison (\ref{g10t}), (\ref{g12a}) and (\ref{e1}) we 
conclude
\bea
\ell_k\left({\bf d}^m\right)={m-1\choose k-1}f_{k-1}\;,\quad 
\ell_1\left({\bf d}^m\right)=1\;.\label{e3}
\eea
\subsection{The 2nd Sylvester wave}\label{s42}
Performing a binomial expansion in (\ref{g9}) and comparing it with (\ref{g10t})
we get
\bea
{\sf L}^{2,1}\left(s,{\bf d}^m\right)=\sum_{r_1,r_2,r_3\geq 0}^{r_1+r_2+r_3=
\omega_2-1}\frac{(\omega_2-1)!}{r_1!\;r_2!\;r_3!}\;l_{r_1}\cdot g_{r_2}\cdot 
s^{r_3},\label{k1}
\eea
where for $l_r$ and $g_r$ the above formula presumes a symbolic exponentiation,
\bea
l_r=\left(\lambda_1+\sum_{i=1}^{\omega_2}B\;d_i\right)^r\!\!,\quad\!g_r=\left(
\gamma_1+\sum_{i=\omega_2+1}^{m}E(0)\;d_i\right)^r\!\!,\!\quad \lambda_k=\sum^{
d_i\in\Xi_2\left({\bf d}^m\right)}_{1\leq i\leq \omega_2}d_i^k,\!\quad\gamma_k=
\sum^{d_i\not\in\Xi_2\left({\bf d}^m\right)}_{\omega_2< i\leq m}d_i^k.\nonumber
\eea
Note that $\sigma_k=\lambda_k+\gamma_k$. Making use of sequences of Bernoulli 
numbers (\ref{e1a}) and values $E_k(0)$ of Euler polynomials,
\bea
E_0(0)=1,\;E_1(0)=-\frac{1}{2},\;E_3(0)=\frac{1}{4},\;E_5(0)=-\frac{1}{2},\;E_7
(0)=\frac{17}{8},\;\ldots,\;\;E_{2k}=0\;,\;k\geq 1,\quad\label{k1a}
\eea
we perform a straightforward calculation in (\ref{k1}) and give the six first 
values $g_r$, $0\leq r\leq 5$,
\bea
&&g_0=1\;,\quad g_1=\frac{\gamma_1}{2}\;,\quad g_2=\frac{1}{2^2}\left(\gamma_1^2-
\gamma_2\right)\;,\quad g_3=\frac{\gamma_1}{2^3}\left(\gamma_1^2-3\gamma_2\right)
\;,\nonumber\\
&&g_4=\frac{1}{2^4}\left(\gamma_1^4-6\gamma_1^2\gamma_2+3\gamma_2^2+2\gamma_4
\right)\;,\quad g_5=\frac{\gamma_1}{2^5}\left(\gamma_1^4-10\gamma_1^2\gamma_2+
15\gamma_2^2+10\gamma_4\right)\;.\nonumber
\eea
For $l_r$ we have to take corresponding $f_r$ given in (\ref{e2}) and replace 
there $\sigma_k$ by $\lambda_k$, i.e., $l_r(\lambda_1,\lambda_2,\ldots)$ 
$\to f_r(\sigma_1,\sigma_2,\ldots)$. By comparison (\ref{g10t}), (\ref{g13}) 
and (\ref{k1}) we conclude
\bea
L^{2,1}_{r+1}\left({\bf d}^m\right)={\omega_2-1\choose r}\sum_{k=0}^r{r\choose 
k}g_k\cdot l_{r-k}\;,\quad\quad0\leq r\leq\omega_2-1,\quad\mbox{i.e.,}
\hspace{1cm}\label{k3}\\
L^{2,1}_1\left({\bf d}^m\right)=1\;,\quad L^{2,1}_2\left({\bf d}^m\right)=
\left(\omega_2-1\right)\frac{\sigma_1}{2}\;,\quad \ldots\;,\quad L^{2,1}_{
\omega_2}\left({\bf d}^m\right)=\sum_{k=0}^{\omega_2-1}{\omega_2-1\choose k}
g_k\cdot l_{\omega_2-1-k}\;.\nonumber
\eea
\subsection{The higher Sylvester waves, $m\geq 3$}\label{s43}
Consider the representation (\ref{g10z}) and rewrite it as follows,
\bea
\begin{array}{r}
{\sf L}^{q,n}\left(s,{\bf d}^m\right)\;=\;{\cal W}_{q,n}\left(s,{\bf d}^m\right)
\xi_q^{sn}\;+\;{\cal W}_{q,-n}\left(s,{\bf d}^m\right)\xi_q^{-sn}\;,\\
i\;{\sf M}^{q,n}\left(s,{\bf d}^m\right)\;=\;{\cal W}_{q,n}\left(s,{\bf d}^m
\right)\xi_q^{sn}\;-\;{\cal W}_{q,-n}\left(s,{\bf d}^m\right)\xi_q^{-sn}\;.
\end{array}\label{y1}
\eea
Substitute a symbolic Rep (\ref{g10e}) into equalities (\ref{y1}) and perform 
their binomial expansions,
\bea
{\sf L}^{q,n}\left(s,{\bf d}^m\right)&=&\sum_{r_1,r_2,r_3\geq 0}^{r_1+r_2+r_3=
\omega_q-1}\frac{(\omega_q-1)!}{r_1!\;r_2!\;r_3!}\left(\Pi^{q,n}_{r_2,+}+
\Pi^{q,n}_{r_2,-}\right)l_{r_1}^{(q)}s^{r_3}\;,\label{y3}\\
i\;{\sf M}^{q,n}\left(s,{\bf d}^m\right)&=&\sum_{r_1,r_2,r_3\geq 0}^{r_1+r_2+
r_3=\omega_q-1}\frac{(\omega_q-1)!}{r_1!\;r_2!\;r_3!}\left(\Pi^{q,n}_{r_2,+}-
\Pi^{q,n}_{r_2,-}\right)l_{r_1}^{(q)}s^{r_3}\;,\nonumber
\eea
where
\bea
\Pi^{q,n}_{r,\pm}=\frac{h_{r,\pm}^{q,n}}{\prod_{i=\omega_q+1}^m\left(1-\xi_q^{
\pm d_i\;n}\right)}\;,\quad\quad h_{r,\pm}^{q,n}=\left(\gamma_1^{(q)}+\sum_{i=
\omega_q+1}^{m}H\left(\xi_q^{\pm d_i\;n}\right)\;d_i\right)^r\;,\label{y4}
\eea
\bea
l_r^{(q)}=\left(\lambda_1^{(q)}+\sum_{i=1}^{\omega_q}B\;d_i\right)^r,\quad
\lambda_k^{(q)}=\sum^{d_i\in\Xi_q\left({\bf d}^m\right)}_{1\leq i\leq \omega_q}
d_i^k,\quad\gamma_k^{(q)}=\sum^{d_i\not\in\Xi_q\left({\bf d}^m\right)}_{
\omega_q< i\leq m}d_i^k,\quad\lambda_k^{(q)}+\gamma_k^{(q)}=\sigma_k.\nonumber
\eea
By comparison (\ref{g10t}), (\ref{g13}), (\ref{g13a}) and (\ref{y3}) we conclude
\bea
L^{q,n}_{r+1}\left({\bf d}^m\right)&=&{\omega_q-1\choose r}\sum_{k=0}^r{r\choose
k}\left(\Pi^{q,n}_{k,+}+\Pi^{q,n}_{k,-}\right)\cdot l_{r-k}^{(q)}\;,\quad\quad
0\leq r\leq\omega_q-1,\hspace{1cm}\label{y6}\\
i\;M^{q,n}_{r+1}\left({\bf d}^m\right)&=&{\omega_q-1\choose r}\sum_{k=0}^r{r
\choose k}\left(\Pi^{q,n}_{k,+}-\Pi^{q,n}_{k,-}\right)\cdot l_{r-k}^{(q)}\;.
\label{y7}
\eea
It is easy to calculate the 1st pair of coefficients
\bea
L^{q,n}_1\left({\bf d}^m\right)=\frac{1+(-1)^{m-\omega_q}\xi_q^{\sigma_1\;n}}
{\prod_{i=\omega_q+1}^m\left(1-\xi_q^{d_i\;n}\right)}\;,\quad\quad
i\;M^{q,n}_1\left({\bf d}^m\right)=\frac{1-(-1)^{m-\omega_q}\xi_q^{\sigma_1\;n}}
{\prod_{i=\omega_q+1}^m\left(1-\xi_q^{d_i\;n}\right)}\;,\label{y8}
\eea
and note that both numbers $L^{q,n}_1\left({\bf d}^m\right)$ and $M^{q,n}_1
\left({\bf d}^m\right)$ cannot vanish simultaneously. Lemma \ref{lem1} provide 
selection rules for these coefficients and explains why some of them disappear 
in (\ref{g5a}).
\vspace{-.2cm}
\begin{lemma}\label{lem1}
\bea
&&\mbox{\rm 1.}\quad\frac{\sigma_1\cdot n}{q}=k,\quad k\in {\mathbb N}\;,  
\label{k13}\\
&&{\sf L}^{q,n}_{1}\left({\bf d}^m\right)=0,\quad\mbox{\rm if}\quad\left\{
\begin{array}{l}2\mid m,\;2\nmid\omega_q,\\2\nmid m,\;2\mid \omega_q,\end{array}
\right.\quad\mbox{\rm and}\quad{\sf M}^{q,n}_{1}\left({\bf d}^m\right)=0,\quad
\mbox{\rm if}\quad\left\{\begin{array}{l}2\mid m,\;2\mid\omega_q,\\2\nmid m,\;
2\nmid \omega_q.\end{array}\right.\nonumber
\eea
\vspace{-.5cm}
\bea
&&\mbox{\rm 2.}\quad\frac{\sigma_1\cdot n}{q}=k+\frac1{2},\quad k\in{\mathbb N}
\;,\label{k14}\\
&&{\sf M}^{q,n}_{1}\left({\bf d}^m\right)=0,\quad\mbox{\rm if}\quad\left\{
\begin{array}{l}2\mid m,\;2\nmid\omega_q,\\2\nmid m,\;2\mid \omega_q,\end{array}
\right.\quad\mbox{\rm and}\quad{\sf L}^{q,n}_{1}\left({\bf d}^m\right)=0,\quad  
\mbox{\rm if}\quad\left\{\begin{array}{l}2\mid m,\;2\mid\omega_q,\\2\nmid m,\;
2\nmid\omega_q.\end{array}\right.\nonumber
\eea
\end{lemma}
\vspace{-.3cm}
\section{Identities for Degrees of the Syzygies}\label{s5}
\vspace{-.3cm}
In this section we bring together two different approaches, theory of restricted
partition and commutative semigroup rings theory, and prove the main Theorems 
\ref{the1} and \ref{the2}. We start with simple identity which comes by 
comparison of (\ref{tr16}) and (\ref{g1}),
\bea
\sum_{s\;\in\;{\sf S}\left({\bf d}^m\right)}z^s=Q\left({\bf d}^m;z\right)\sum_{
s\;\in\;{\sf S}\left({\bf d}^m\right)}W\left({\bf d}^m;s\right)z^s\;.\label{z3}
\eea
Substituting (\ref{bet05}) and (\ref{bet1}) into (\ref{z3}) we get
\bea
\sum_{s\;\in\;{\sf S}\left({\bf d}^m\right)}z^s&=&\sum_{s\;\in\;{\sf S}\left(
{\bf d}^m\right)}W\left({\bf d}^m;s\right)z^s-\sum_{j=1}^{\beta_1\left({\bf d}^m
\right)}\sum_{s\;\in\;{\sf S}\left({\bf d}^m\right)}W\left({\bf d}^m;s\right)
z^{C_{j,1}+s}+\ldots+\nonumber\\
&&(-1)^{m-1}\sum_{j=1}^{\beta_{m-1}\left({\bf d}^m\right)}\sum_{s\;\in\;{\sf S}
\left({\bf d}^m\right)}W\left({\bf d}^m;s\right)z^{C_{j,m-1}+s}\;.\label{z5}
\eea
Rewrite an equality (\ref{z5}) as follows
\bea
\sum_{s\;\in\;{\sf S}\left({\bf d}^m\right)}z^s&=&\sum_{s\;\in\;{\sf S}\left(
{\bf d}^m\right)}W\left({\bf d}^m;s\right)z^s-\sum_{j=1}^{\beta_1\left({\bf d}^m
\right)}\sum_{s\;\in\;{\sf S}\left({\bf d}^m\right)}W\left({\bf d}^m;s-C_{j,1}
\right)z^s+\ldots+\nonumber\\
&&(-1)^{m-1}\sum_{j=1}^{\beta_{m-1}\left({\bf d}^m\right)}\sum_{s\;\in\;{\sf S}
\left({\bf d}^m\right)}W\left({\bf d}^m;s-C_{j,m-1}\right)z^s\;,\label{z6}
\eea
and equate the corresponding contributions coming from monomial terms $z^s$ in 
the left hand side (l.h.s.) and right hand side (r.h.s.) of (\ref{z6}). This 
gives a quasipolynomial equality,
\bea
W\left({\bf d}^m;s\right)-\sum_{j=1}^{\beta_1\left({\bf d}^m\right)}W\left({\bf 
d}^m;s-C_{j,1}\right)+\ldots +(-1)^{m-1}\sum_{j=1}^{\beta_{m-1}\left({\bf d}^m
\right)}W\left({\bf d}^m;s-C_{j,m-1}\right)=1\;.\label{z7}
\eea
Now substitute into (\ref{z7}) the Sylvester expansions (\ref{g6}), (\ref{g10r})
and make use of the linear independence of the partial Sylvester waves $W_q\!
\left(s,{\bf d}^m\right)$ and $\overline{{\cal W}}_{q,n}\!\left(s,{\bf d}^m
\right)$. 

This gives rise to a set of quasipolynomial equalities,
\bea
W_1\!\left(s,{\bf d}^m\right)-\sum_{j=1}^{\beta_1\left({\bf d}^m\right)}W_1\!
\left(s-C_{j,1},{\bf d}^m\right)+\ldots +(-1)^{m-1}\sum_{j=1}^{\beta_{m-1}\left(
{\bf d}^m\right)}W_1\!\left(s-C_{j,m-1},{\bf d}^m\right)=1,\quad\label{z8}
\eea
and for $q\geq 2$, $1\leq n<q/2$ such that $\gcd(n,q)=1$,
\bea
\overline{{\cal W}}_{q,n}\!\left(s,{\bf d}^m\right)-\!\sum_{j=1}^{\beta_1\left(
{\bf d}^m\right)}\!\!\overline{{\cal W}}_{q,n}\!\left(s-C_{j,1},{\bf d}^m
\right)+\ldots +(-1)^{m-1}\!\!\sum_{j=1}^{\beta_{m-1}\left({\bf d}^m\right)}\!
\!\overline{{\cal W}}_{q,n}\left(s-C_{j,m-1},{\bf d}^m\right)=0.\quad\label{z9}
\eea
Equalities (\ref{z8}) and (\ref{z9}) are a source of new relationships between
degrees $C_{j,i}$ of the syzygies of different kinds. What is more important 
that for this purpose we have to know only the most basic properties of partial 
Sylvester waves $W_q\!\left(s,{\bf d}^m\right)$ such as their quasipolynomial 
Reps (\ref{g10r}), (\ref{g10y}) -- (\ref{g13a}) and the coefficients $\ell_1
\left({\bf d}^m\right)$, $L_1^{q,n}\left({\bf d}^m\right)$ and $M_1^{q,n}\left(
{\bf d}^m\right)$ at the leading terms (\ref{e3}), (\ref{k3}) and (\ref{y8}) but
not the whole set of quasipolynomial coefficients in (\ref{g12a}) -- 
(\ref{g13a}).
\subsection{Polynomial Identities Associated with the 1st Sylvester Wave.\\
The Proof of Theorem \ref{the1}}\label{s51}
In this section we prove Theorem \ref{the1} which was stated in section 
\ref{s21}. Consider equality (\ref{z8}) and substitute Rep (\ref{g10t}) for 
$W_1\!\left(s,{\bf d}^m\right)$ into (\ref{z8}),
\bea
{\sf L}^1\left(s,{\bf d}^m\right)-\sum_{j=1}^{\beta_1\left({\bf d}^m\right)}
{\sf L}^1\!\left(s-C_{j,1},{\bf d}^m\right)+\ldots +(-1)^{m-1}\sum_{j=1}^{
\beta_{m-1}\left({\bf d}^m\right)}{\sf L}^1\!\left(s-C_{j,m-1},{\bf d}^m\right)=
(m-1)!\;\pi_m.\nonumber
\eea
Inserting Rep (\ref{g12a}) into the last identity we get
\bea
&&\ell_1\left({\bf d}^m\right)\left[s^{m-1}-\sum_{j=1}^{\beta_1\left({\bf d}^m
\right)}\left(s-C_{j,1}\right)^{m-1}+\ldots+(-1)^{m-1}\sum_{j=1}^{\beta_{m-1}
\left({\bf d}^m\right)}\left(s-C_{j,m-1}\right)^{m-1}\right]+\nonumber\\
&&\ell_2\left({\bf d}^m\right)\left[s^{m-2}-\sum_{j=1}^{\beta_1\left({\bf d}^m
\right)}\left(s-C_{j,1}\right)^{m-2}+\ldots+(-1)^{m-1}\sum_{j=1}^{\beta_{m-2}
\left({\bf d}^m\right)}\left(s-C_{j,m-1}\right)^{m-2}\right]+\ldots +\nonumber\\
&&\ell_{m-1}\left({\bf d}^m\right)\left[s-\sum_{j=1}^{\beta_1\left({\bf d}^m
\right)}\left(s-C_{j,1}\right)+\ldots+(-1)^{m-1}\sum_{j=1}^{\beta_{m-1}\left(
{\bf d}^m\right)}\left(s-C_{j,m-1}\right)\right]+\nonumber\\
&&\ell_m\left({\bf d}^m\right)\left[1-\sum_{j=1}^{\beta_1\left({\bf d}^m\right)}
1+\ldots+(-1)^{m-1}\sum_{j=1}^{\beta_{m-1}\left({\bf d}^m\right)}1\right]=
(m-1)!\;\pi_m\;.\label{z10}
\eea
Introduce the following sums and for short denote them as follows,
\bea
{\bf A}_0\left({\bf d}^m\right)&=&1-\beta_1\left({\bf d}^m\right)+\beta_2\left(
{\bf d}^m\right)-\ldots+(-1)^{m-1}\beta_{m-1}\left({\bf d}^m\right)\;,
\nonumber\\
{\bf A}_k\left({\bf d}^m\right)&=&\sum_{j=1}^{\beta_1\left({\bf d}^m\right)}
C_{j,1}^k-\sum_{j=1}^{\beta_2\left({\bf d}^m\right)}C_{j,2}^k+\ldots+(-1)^m
\sum_{j=1}^{\beta_{m-1}\left({\bf d}^m\right)}C_{j,m-1}^k\;.\label{z11}
\eea
Keeping in mind (\ref{z11}) rewrite equality (\ref{z10}) and get
\bea
&&\left[\ell_1\left({\bf d}^m\right)s^{m-1}+\ell_2\left({\bf d}^m\right)s^{m-2}
+\ldots+\ell_m\left({\bf d}^m\right)\right]{\bf A}_0\left({\bf d}^m
\right)+\quad\quad\label{z11a}\\
&&\left[\ell_1\left({\bf d}^m\right){m-1\choose 1}s^{m-2}+\ell_2\left({\bf d}^m
\right){m-2\choose 1}s^{m-3}+\ldots+\ell_{m-1}\left({\bf d}^m\right)
{1\choose 1}\right]{\bf A}_1\left({\bf d}^m\right)-\nonumber\\
&&\left[\ell_1\left({\bf d}^m\right){m-1\choose 2}s^{m-3}+\ell_2\left({\bf d}^m
\right){m-2\choose 2}s^{m-4}+\ldots+\ell_{m-2}\left({\bf d}^m\right)
{2\choose 2}\right]{\bf A}_2\left({\bf d}^m\right)+\ldots-\nonumber\\
&&(-1)^m\left[\ell_1\left({\bf d}^m\right){m-1\choose m-2}s+\ell_2\left(
{\bf d}^m\right)\right]{\bf A}_{m-2}\left({\bf d}^m\right)+\nonumber\\
&&(-1)^m\;\ell_1\left({\bf d}^m\right){\bf A}_{m-1}\left({\bf d}^m\right)=
(m-1)!\pi_m\;.\nonumber
\eea
Equating the corresponding contributions coming from the power terms $s^a$,
$0\leq a<m$, in the l.h.s. and the r.h.s. of (\ref{z11a}) and keeping in 
mind $\ell_1\left({\bf d}^m\right)=1$ (see (\ref{e3})), we get finally,
\begin{eqnarray}
{\bf A}_k\left({\bf d}^m\right)=0\;,\quad k=0,1,\ldots,m-2\;,\quad\quad
{\bf A}_{m-1}\left({\bf d}^m\right)=(-1)^m(m-1)!\;\pi_m\;.\label{z13}
\eea
The 1st identity ${\bf A}_0\left({\bf d}^m\right)=0$ is already known in 
(\ref{bet2}). The rest of identities prove Theorem \ref{the1}.
\vspace{-1cm}
\subsection{Quasipolynomial Identities Associated with the 2nd Sylvester Wave}
\label{s52}
In this section we study an intermediate case $q=2$ of the master equality 
(\ref{z9}) which is technically slightly more difficult than the previous 
equality (\ref{z8}). Consider (\ref{z9}) and substitute Rep (\ref{g10t}) for 
$W_2\!\left(s,{\bf d}^m\right)=\overline{{\cal W}}_{2,1}\left(s,{\bf d}^m
\right)$ into (\ref{z9}),
\bea
{\sf L}^{2,1}\left(s,{\bf d}^m\right)\cos\pi s-\sum_{j=1}^{\beta_1\left({\bf 
d}^m\right)}{\sf L}^{2,1}\!\left(s-C_{j,1},{\bf d}^m\right)\cos\pi\left(s-
C_{j,1}\right)+\ldots+\hspace{2cm}\nonumber\\
(-1)^{m-1}\sum_{j=1}^{\beta_{m-1}\left({\bf d}^m\right)}{\sf L}^{2,1}\!\left(s
-C_{j,m-1},{\bf d}^m\right)\cos\pi\left(s-C_{j,m-1}\right)=0\;.\label{z14}
\eea
Substituting Rep (\ref{g13}) into the last equality we get
\bea
\left(L_1^{2,1}\left({\bf d}^m\right)s^{\omega_2-1}+L_2^{2,1}\left({\bf d}^m
\right)s^{\omega_2-2}+\ldots+L_{\omega_2-1}^{2,1}\left({\bf d}^m\right)s+
L_{\omega_2}^{2,1}\left({\bf d}^m\right)\right)\cos\pi s-\hspace{1.5cm}
\nonumber\\
\sum_{j=1}^{\beta_1\left({\bf d}^m\right)}\left(L_1^{2,1}\left({\bf d}^m\right)
\left(s-C_{j,1}\right)^{\omega_2-1}+\ldots+L_{\omega_2}^{2,1}\left({\bf d}^m
\right)\right)\cos\pi\left(s-C_{j,1}\right)+\hspace{1cm}\nonumber\\
\ldots+(-1)^{m-1}\sum_{j=1}^{\beta_{m-1}\left({\bf d}^m\right)}\left(
L_1^{2,1}\left({\bf d}^m\right)\left(s-C_{j,m-1}\right)^{\omega_2-1}+\ldots+
L_{\omega_2}^{2,1}\left({\bf d}^m\right)\right)\cos\pi\left(s-C_{j,m-1}\right)
=0\;.\nonumber
\eea
Keeping in mind two identities $\sin\pi s=0$ and $\cos\pi s\neq0$, 
$s\in{\mathbb N}$, we write,
\bea
\left(L_1^{2,1}\left({\bf d}^m\right)s^{\omega_2-1}+L_2^{2,1}\left({\bf d}^m
\right)s^{\omega_2-2}+\ldots+L_{\omega_2-1}^{2,1}\left({\bf d}^m\right)s+
L_{\omega_2}^{2,1}\left({\bf d}^m\right)\right)-\hspace{1.5cm}\label{z12}\\
\sum_{j=1}^{\beta_1\left({\bf d}^m\right)}\left(L_1^{2,1}\left({\bf d}^m\right)
\left(s-C_{j,1}\right)^{\omega_2-1}+\ldots+L_{\omega_2}^{2,1}\left({\bf d}^m
\right)\right)\cos\left(\pi C_{j,1}\right)+\ldots+\hspace{1cm}\nonumber\\
(-1)^{m-1}\sum_{j=1}^{\beta_{m-1}\left({\bf d}^m\right)}\left(L_1^{2,1}\left(
{\bf d}^m\right)\left(s-C_{j,m-1}\right)^{\omega_2-1}+\ldots+L_{\omega_2}^{
2,1}\left({\bf d}^m\right)\right)\cos\left(\pi C_{j,m-1}\right)=0.\nonumber
\eea
Rearrange the last equality (\ref{z12}) as follows,
\bea
&&L_1^{2,1}\left({\bf d}^m\right)\left[s^{\omega_2-1}-\sum_{j=1}^{\beta_1\left(
{\bf d}^m\right)}\left(s-C_{j,1}\right)^{\omega_2-1}\cos\left(\pi C_{j,1}
\right)+\ldots+\right.\label{z15}\\
&&\left.\hspace{3cm}(-1)^{m-1}\sum_{j=1}^{\beta_{m-1}\left({\bf d}^m\right)}
\left(s-C_{j,m-1}\right)^{\omega_2-1}\cos\left(\pi C_{j,m-1}\right)\right]+
\nonumber\\
&&L_2^{2,1}\left({\bf d}^m\right)\left[s^{\omega_2-2}-\sum_{j=1}^{\beta_1\left(
{\bf d}^m\right)}\left(s-C_{j,1}\right)^{\omega_2-2}\cos\left(\pi C_{j,1}\right)
+\ldots+\right.\nonumber\\
&&\left.\hspace{3cm}(-1)^{m-1}\sum_{j=1}^{\beta_{m-2}\left({\bf d}^m\right)}
\left(s-C_{j,m-1}\right)^{\omega_2-2}\cos\left(\pi C_{j,m-1}\right)\right]+
\ldots +\nonumber\\
&&L_{\omega_2}^{2,1}\left({\bf d}^m\right)\left[1-\sum_{j=1}^{\beta_1\left({\bf 
d}^m\right)}\cos\left(\pi C_{j,1}\right)+\ldots+(-1)^{m-1}\sum_{j=1}^{\beta_{
m-1}\left({\bf d}^m\right)}\cos\left(\pi C_{j,m-1}\right)\right]=0\;.
\nonumber
\eea
and denote two following sums,
\bea
{\bf B}_0\left({\bf d}^m\right)&=&1-\sum_{j=1}^{\beta_1\left({\bf d}^m\right)}
\cos\left(\pi C_{j,1}\right)+\ldots+(-1)^{m-1}\sum_{j=1}^{\beta_{m-1}
\left({\bf d}^m\right)}\cos\left(\pi C_{j,m-1}\right)\;,\label{z16}\\
{\bf B}_k\left({\bf d}^m\right)&=&\sum_{j=1}^{\beta_1\left({\bf d}^m\right)}
C_{j,1}^k\cos\left(\pi C_{j,1}\right)-\ldots+(-1)^m\sum_{j=1}^{\beta_{m-1}
\left({\bf d}^m\right)}C_{j,m-1}^k\cos\left(\pi C_{j,m-1}\right).\nonumber
\eea
Keeping in mind (\ref{z16}) rewrite equality (\ref{z15}) and get
\bea
&&\left[L^{2,1}_1\left({\bf d}^m\right)s^{\omega_2-1}+L^{2,1}_2\left({\bf d}^m
\right)s^{\omega_2-2}+\ldots+L^{2,1}_{\omega_2}\left({\bf d}^m\right)\right]
{\bf B}_0\left({\bf d}^m\right)+\label{z17}\\
&&\left[L^{2,1}_1\left({\bf d}^m\right){\omega_2-1\choose 1}s^{\omega_2-2}+
L^{2,1}_2\left({\bf d}^m\right){\omega_2-2\choose 1}s^{\omega_2-3}+\ldots+
L^{2,1}_{\omega_2-1}\left({\bf d}^m\right)\right]{\bf B}_1\left({\bf d}^m
\right)-\nonumber\\
&&\left[L^{2,1}_1\left({\bf d}^m\right){\omega_2-1\choose 2}s^{\omega_2-3}+
L^{2,1}_2\left({\bf d}^m\right){\omega_2-2\choose 2}s^{\omega_2-4}+\ldots+
L^{2,1}_{\omega_2-2}\left({\bf d}^m\right)\right]{\bf B}_2\left({\bf d}^m
\right)+\ldots-\nonumber\\
&&(-1)^{\omega_2}\left[L^{2,1}_1\left({\bf d}^m\right){\omega_2-1\choose
\omega_2-2}s+L^{2,1}_2\left({\bf d}^m\right)\right]{\bf B}_{\omega_2-2}
\left({\bf d}^m\right)+\nonumber\\
&&(-1)^{\omega_2}L^{2,1}_1\left({\bf d}^m\right){\bf B}_{\omega_2-1}
\left({\bf d}^m\right)=0\;.\nonumber
\eea
Equating the corresponding contributions coming from the power terms $s^a$, $0
\leq a<\omega_2$, in the l.h.s. and the r.h.s. of (\ref{z17}) and keeping in
mind $L^{2,1}_1\left({\bf d}^m\right)=1$ (see (\ref{k3})), we get finally,
\bea
{\bf B}_k\left({\bf d}^m\right)=0\;,\quad k=0,1,\ldots,\omega_2-1\;.\label{z18}
\eea
\subsection{Quasipolynomial Identities Associated with the Higher Sylvester 
Waves.\\The Proof of Theorem \ref{the2}}\label{s53}
In this section we study the general case of equality (\ref{z9}) and start to 
prove Theorem \ref{the2}, we finish its proof in section \ref{s533}. Consider
(\ref{z9}) and substitute Reps (\ref{g10r}) and (\ref{g10y}) into (\ref{z9}),
\bea
&&{\sf L}^{q,n}\left(s,{\bf d}^m\right)\cos\frac{2\pi n}{q}s+{\sf M}^{q,n}
\left(s,{\bf d}^m\right)\sin\frac{2\pi n}{q}s-\nonumber\\
&&\sum_{j=1}^{\beta_1\left({\bf d}^m\right)}\left({\sf L}^{q,n}\!\left(s-
C_{j,1},{\bf d}^m\right)\cos\frac{2\pi n}{q}\left(s-C_{j,1}\right)+
{\sf M}^{q,n}\!\left(s-C_{j,1},{\bf d}^m\right)\sin\frac{2\pi n}{q}\left(s-
C_{j,1}\right)\right)+\nonumber\\
&&\sum_{j=1}^{\beta_2\left({\bf d}^m\right)}\left({\sf L}^{q,n}\!\left(s-
C_{j,2},{\bf d}^m\right)\cos\frac{2\pi n}{q}\left(s-C_{j,2}\right)+
{\sf M}^{q,n}\!\left(s-C_{j,2},{\bf d}^m\right)\sin\frac{2\pi n}{q}\left(s-
C_{j,2}\right)\right)-\nonumber\\
&&\hspace{1cm}\ldots+(-1)^{m-1}\sum_{j=1}^{\beta_{m-1}\left({\bf d}^m\right)}
\left({\sf L}^{q,n}\!\left(s-C_{j,m-1},{\bf d}^m\right)\cos\frac{2\pi n}{q}
\left(s-C_{j,m-1}\right)+\right.\nonumber\\
&&\left.\hspace{4.5cm}{\sf M}^{q,n}\!\left(s-C_{j,1},{\bf d}^m\right)\sin
\frac{2\pi n}{q}\left(s-C_{j,1},{\bf d}^m\right)\right)=0\;.\nonumber
\eea
For  $q\geq 2$, $1\leq n<q/2$, $\gcd(n,q)=1$, we represent the last equality in 
the form
\bea
{\mathbb A}^{q,n}\left(s,{\bf d}^m\right)\cos\frac{2\pi n}{q}s+  
{\mathbb B}^{q,n}\left(s,{\bf d}^m\right)\sin\frac{2\pi n}{q}s=0\;,\label{z21a}
\eea
where ${\mathbb A}^{q,n}\left(s,{\bf d}^m\right)$ and ${\mathbb B}^{q,n}\left(
s,{\bf d}^m\right)$ are two polynomials in $s$ which can be calculated,
\bea 
{\mathbb A}^{q,n}\left(s,{\bf d}^m\right)=
{\sf L}^{q,n}\left(s,{\bf d}^m\right)-\sum_{j=1}^{\beta_1\left({\bf d}^m\right)}
{\sf L}^{q,n}\!\left(s-C_{j,1},{\bf d}^m\right)\cos\frac{2\pi n}{q}C_{j,1}+
\hspace{2cm}\nonumber\\
\sum_{j=1}^{\beta_2\left({\bf d}^m\right)}\!\!\!\!{\sf L}^{q,n}\!\left(s-
C_{j,2},{\bf d}^m\right)\cos\frac{2\pi n}{q}C_{j,2}-\ldots +
(-1)^{m-1}\!\!\!\sum_{j=1}^{\beta_{m-1}\left({\bf d}^m\right)}\!\!\!\!{\sf L}^{
q,n}\!\left(s-C_{j,m-1},{\bf d}^m\right)\cos\frac{2\pi n}{q}C_{j,m-1}+
\nonumber\\
\sum_{j=1}^{\beta_1\left({\bf d}^m\right)}\!\!{\sf M}^{q,n}\!\left(s-
C_{j,1},{\bf d}^m\right)\sin\frac{2\pi n}{q}C_{j,1}-\ldots-(-1)^{m-1}\!\!
\sum_{j=1}^{\beta_{m-1}\left({\bf d}^m\right)}\!\!\!\!{\sf M}^{q,n}
\left(s-C_{j,m-1},{\bf d}^m\right)\sin\frac{2\pi n}{q}C_{j,m-1}\nonumber
\eea
\bea
{\mathbb B}^{q,n}\left(s,{\bf d}^m\right)=
{\sf M}^{q,n}\left(s,{\bf d}^m\right)-\sum_{j=1}^{\beta_1\left({\bf d}^m\right)}
{\sf M}^{q,n}\!\left(s-C_{j,1},{\bf d}^m\right)\cos\frac{2\pi n}{q}C_{j,1}+
\hspace{2cm}\nonumber\\
\sum_{j=1}^{\beta_2\left({\bf d}^m\right)}\!\!\!\!{\sf M}^{q,n}\!\left(s-
C_{j,2},{\bf d}^m\right)\cos\frac{2\pi n}{q}C_{j,2}-\ldots +\!(-1)^{m-1}
\!\!\!\!\sum_{j=1}^{\beta_{m-1}\left({\bf d}^m\right)}\!\!\!\!\!{\sf M}^{
q,n}\!\left(s-C_{j,m-1},{\bf d}^m\right)\cos\frac{2\pi n}{q}C_{j,m-1}-
\nonumber\\
\sum_{j=1}^{\beta_1\left({\bf d}^m\right)}\!\!{\sf L}^{q,n}\!\left(s-
C_{j,1},{\bf d}^m\right)\sin\frac{2\pi n}{q}C_{j,1}+\ldots+(-1)^{m-1}\!\!
\sum_{j=1}^{\beta_{m-1}\left({\bf d}^m\right)}\!\!\!\!{\sf L}^{q,n}
\left(s-C_{j,m-1},{\bf d}^m\right)\sin\frac{2\pi n}{q}C_{j,m-1}\nonumber
\eea
The quasipolynomial in the l.h.s. of (\ref{z21a}) vanishes identically iff all 
coefficients at every harmonics $\cos\frac{2\pi n}{q}s$ and $\sin\frac{2\pi n}
{q}s$ do vanish. This is why (\ref{z21a}) is equivalent two independent Eqs,
\bea
{\mathbb A}^{q,n}\left(s,{\bf d}^m\right)=0\;,\quad\quad
{\mathbb B}^{q,n}\left(s,{\bf d}^m\right)=0\;,\label{z22a}
\eea
which have to be satisfied identically but solved separately. 
\subsubsection{${\mathbb A}^{q,n}\left(s,{\bf d}^m\right)=0$}\label{s531}
In this section we find the 1st system of linear equations (\ref{z26}) for the 
coefficients $L_k^{q,n}\left({\bf d}^m\right)$ and $M_k^{q,n}\left({\bf d}^m
\right)$, $1\leq k\leq \omega_q$. First, in order to avoid the lengthy formulas 
coming from (\ref{z22a}) we introduce another pair of quasipolynomials equipped 
with subscript $p$, $0\leq p<\omega_q$,
\bea
{\mathbb C}^{q,n}_p\left(s,{\bf d}^m\right)&=&\sum_{j=1}^{\beta_1\left({\bf d}^m
\right)}\left(s-C_{j,1}\right)^p\cos\frac{2\pi n}{q}C_{j,1}-\sum_{j=1}^{\beta_2
\left({\bf d}^m\right)}\left(s-C_{j,2}\right)^p\cos\frac{2\pi n}{q}C_{j,2}+
\ldots-\nonumber\\
&&\hspace{2cm}(-1)^{m-1}\sum_{j=1}^{\beta_{m-1}\left({\bf d}^m\right)}
\left(s-C_{j,m-1}\right)^p\cos\frac{2\pi n}{q}C_{j,m-1}\;,\label{h1}
\eea
\bea
{\mathbb S}^{q,n}_p\left(s,{\bf d}^m\right)&=&\sum_{j=1}^{\beta_1\left({\bf d}^m
\right)}\left(s-C_{j,1}\right)^p\sin\frac{2\pi n}{q}C_{j,1}-\sum_{j=1}^{\beta_2
\left({\bf d}^m\right)}\left(s-C_{j,2}\right)^p\sin\frac{2\pi n}{q}C_{j,2}+
\ldots-\nonumber\\
&&\hspace{2cm}(-1)^{m-1}\sum_{j=1}^{\beta_{m-1}\left({\bf d}^m\right)}
\left(s-C_{j,m-1}\right)^p\sin\frac{2\pi n}{q}C_{j,m-1}\;.\label{h2}
\eea
Substitute now Reps (\ref{g13}) and (\ref{g13a}) into expression for ${\mathbb 
A}^{q,n}\left(s,{\bf d}^m\right)$ given in the previous section and make use of 
(\ref{h1}) and (\ref{h2}). Thus, we present the 1st Eq. in (\ref{z22a}) as 
follows,
\bea
&&L_1^{q,n}\left({\bf d}^m\right)\left[s^{\omega_q-1}-{\mathbb C}^{q,n}_{
\omega_q-1}\left(s,{\bf d}^m\right)\right]+L_2^{q,n}\left({\bf d}^m\right)
\left[s^{\omega_q-2}-{\mathbb C}^{q,n}_{\omega_q-2}\left(s,{\bf d}^m\right)
\right]+\ldots+\quad\label{h3}\\
&&L_{\omega_q-1}^{q,n}\left({\bf d}^m\right)\left[s-{\mathbb C}^{q,n}_1\left(s,
{\bf d}^m\right)\right]+L_{\omega_q}^{q,n}\left({\bf d}^m\right)\left[1-
{\mathbb C}^{q,n}_0\left(s,{\bf d}^m\right)\right]+\nonumber\\
&&M_1^{q,n}\left({\bf d}^m\right){\mathbb S}^{q,n}_{\omega_q-1}\left(s,{\bf d}^m
\right)+M_2^{q,n}\left({\bf d}^m\right){\mathbb S}^{q,n}_{\omega_q-2}\left(
s,{\bf d}^m\right)+\ldots+M_{\omega_q}^{q,n}\left({\bf d}^m\right){\mathbb S}^{
q,n}_0\left(s,{\bf d}^m\right)=0\;.\hspace{1cm}\nonumber
\eea
After binomial expansion of the terms $\left(s-C_{j,i}\right)^p$ in (\ref{h1}) 
and (\ref{h2}) and subsequent substitution into (\ref{h3}) we recast the power 
terms in $s$ and get the final equality
\bea
&&\left[L^{q,n}_1\left({\bf d}^m\right)s^{\omega_q-1}+L^{q,n}_2\left({\bf d}^m
\right)s^{\omega_q-2}+\ldots+L^{q,n}_{\omega_q}\left({\bf d}^m\right)\right]  
{\bf G}_0\left({\bf d}^m\right)+\label{z25}\\
&&\left[L^{q,n}_1\left({\bf d}^m\right){\omega_q-1\choose 1}s^{\omega_q-2}+
L^{q,n}_2\left({\bf d}^m\right){\omega_q-2\choose 1}s^{\omega_q-3}+\ldots+
L^{q,n}_{\omega_q-1}\left({\bf d}^m\right)\right]{\bf G}_1\left({\bf d}^m
\right)-\nonumber\\
&&\left[L^{q,n}_1\left({\bf d}^m\right){\omega_q-1\choose 2}s^{\omega_q-3}+\!
L^{q,n}_2\left({\bf d}^m\right){\omega_q-2\choose 2}s^{\omega_q-4}+\!\ldots+\!
L^{q,n}_{\omega_q-2}\left({\bf d}^m\right)\right]{\bf G}_2\left({\bf d}^m
\right)+\!\ldots\!-\nonumber\\
&&(-1)^{\omega_q}\left[L^{q,n}_1\left({\bf d}^m\right){\omega_q-1
\choose\omega_q-2}s+L^{q,n}_2\left({\bf d}^m\right)\right]{\bf G}_{\omega_q-2}
\left({\bf d}^m\right)+
(-1)^{\omega_q}\;L^{q,n}_1\left({\bf d}^m\right){\bf G}_{\omega_q-1}\left(
{\bf d}^m\right)+\nonumber\\
&&\left[M^{q,n}_1\left({\bf d}^m\right)s^{\omega_q-1}+M^{q,n}_2\left({\bf d}^m
\right)s^{\omega_q-2}+\ldots+M^{q,n}_{\omega_q}\left({\bf d}^m\right)\right]
{\bf D}_0\left({\bf d}^m\right)-\nonumber\\
&&\left[M^{q,n}_1\left({\bf d}^m\right){\omega_q-1\choose 1}s^{\omega_q-2}+
M^{q,n}_2\left({\bf d}^m\right){\omega_q-2\choose 1}s^{\omega_q-3}+\ldots+
M^{q,n}_{\omega_q-1}\left({\bf d}^m\right)\right]{\bf D}_1\left({\bf d}^m
\right)+\nonumber\\
&&\left[M^{q,n}_1\left({\bf d}^m\right){\omega_q-1\choose 2}s^{\omega_q-3}+\!
M^{q,n}_2\left({\bf d}^m\right){\omega_q-2\choose 2}s^{\omega_q-4}+\!\ldots\!+
M^{q,n}_{\omega_q-2}\left({\bf d}^m\right)\right]{\bf D}_2\left({\bf d}^m
\right)-\!\ldots\!+\nonumber\\
&&(-1)^{\omega_q}\left[M^{q,n}_1\left({\bf d}^m\right){\omega_q-1\choose
\omega_q-2}s+\!M^{q,n}_2\left({\bf d}^m\right)\right]{\bf D}_{\omega_q-2}
\left({\bf d}^m\right)-
(-1)^{\omega_q}M^{q,n}_1\left({\bf d}^m\right){\bf D}_{\omega_q-1}\left(
{\bf d}^m\right)=0,\nonumber
\eea
where we denote for short the following sums,
\bea
&&{\bf G}_0\left({\bf d}^m\right)=1-\sum_{j=1}^{\beta_1\left({\bf d}^m\right)}
\cos\frac{2\pi n}{q}C_{j,1}+\ldots+(-1)^{m-1}\sum_{j=1}^{\beta_{m-1}
\left({\bf d}^m\right)}\cos\frac{2\pi n}{q}C_{j,m-1}\;,\label{z24}\\
&&{\bf G}_k\left({\bf d}^m\right)=\sum_{j=1}^{\beta_1\left({\bf d}^m\right)}
C_{j,1}^k\cos\frac{2\pi n}{q}C_{j,1}-\ldots+(-1)^m\sum_{j=1}^{\beta_{m-1}
\left({\bf d}^m\right)}C_{j,m-1}^k\cos\frac{2\pi n}{q}C_{j,m-1}\;,\nonumber\\
&&{\bf D}_k\left({\bf d}^m\right)=\sum_{j=1}^{\beta_1\left({\bf d}^m\right)}
C_{j,1}^k\sin\frac{2\pi n}{q}C_{j,1}-\ldots+(-1)^m\sum_{j=1}^{\beta_{m-1}
\left({\bf d}^m\right)}C_{j,m-1}^k\sin\frac{2\pi n}{q}C_{j,m-1}\;.\nonumber
\eea
Equating the corresponding contributions coming from the power terms $s^a$, 
$a=0,\ldots,\omega_q-1$, in the l.h.s. and the r.h.s. of (\ref{z25}) we get 
a system of $\omega_q$ linear equations for $L_k^{q,n}\left({\bf d}^m\right)$ 
and $M_k^{q,n}\left({\bf d}^m\right)$. For short we present here the three 
first equations
\bea
&&a=\omega_q-1\hspace{2cm}
L^{q,n}_1\left({\bf d}^m\right){\bf G}_0\left({\bf d}^m\right)+
M^{q,n}_1\left({\bf d}^m\right){\bf D}_0 \left({\bf d}^m\right)=0\;,
\hspace{1cm}\label{z27}\\ \nonumber\\
&&a=\omega_q-2\hspace{2cm}
L^{q,n}_2\left({\bf d}^m\right){\bf G}_0\left({\bf d}^m\right)+
M^{q,n}_2\left({\bf d}^m\right){\bf D}_0\left({\bf d}^m\right)+\label{z27a}\\
&&\hspace{4cm}{\omega_q-1\choose 1}\left(L^{q,n}_1\left({\bf d}^m\right)
{\bf G}_1\left({\bf d}^m\right)-M^{q,n}_1\left({\bf d}^m\right){\bf D}_1
\left({\bf d}^m\right)\right)=0\;,\hspace{1cm}\nonumber\\\nonumber\\
&&a=\omega_q-3\hspace{2cm}
L^{q,n}_3\left({\bf d}^m\right){\bf G}_0\left({\bf d}^m\right)+
M^{q,n}_3\left({\bf d}^m\right){\bf D}_0\left({\bf d}^m\right)+\nonumber\\
&&\hspace{4cm}{\omega_q-2\choose 1}\left(L^{q,n}_2\left({\bf d}^m\right)
{\bf G}_1\left({\bf d}^m\right)-M^{q,n}_2\left({\bf d}^m\right)  
{\bf D}_1\left({\bf d}^m\right)\right)-\nonumber\\
&&\hspace{4cm}{\omega_q-1\choose 2}\left(L^{q,n}_1\left({\bf d}^m\right)
{\bf G}_2\left({\bf d}^m\right)-M^{q,n}_1\left({\bf d}^m\right)   
{\bf D}_2\left({\bf d}^m\right)\right)=0\;,\hspace{1cm}\nonumber
\eea
and the last one
\bea
&&a=0\hspace{1.5cm}
L^{q,n}_{\omega_q}\left({\bf d}^m\right){\bf G}_0\left({\bf d}^m\right)
+M^{q,n}_{\omega_q}\left({\bf d}^m\right){\bf D}_0\left({\bf d}^m\right)+
\nonumber\\
&&\hspace{2.4cm}
L^{q,n}_{\omega_q-1}\left({\bf d}^m\right){\bf G}_1\left({\bf d}^m\right)-
M^{q,n}_{\omega_q-1}\left({\bf d}^m\right){\bf D}_1\left({\bf d}^m\right)-
\ldots-\nonumber\\
&&\hspace{2.5cm}
(-1)^{\omega_q}\left(L^{q,n}_2\left({\bf d}^m\right){\bf G}_{\omega_q-2}
\left({\bf d}^m\right)-M^{q,n}_2\left({\bf d}^m\right){\bf D}_{\omega_q-2}
\left({\bf d}^m\right)\right)+\nonumber\\
&&\hspace{2.5cm}(-1)^{\omega_q}\left(L^{q,n}_1\left({\bf d}^m\right)
{\bf G}_{\omega_q-1}\left({\bf d}^m\right)-M^{q,n}_1\left({\bf d}^m\right)
{\bf D}_{\omega_q-1}\left({\bf d}^m\right)\right)=0\;.\nonumber
\eea
The whole set of these equations can be written in more compact form
\bea
&&L^{q,n}_k\left({\bf d}^m\right){\bf G}_0\left({\bf d}^m\right)+
M^{q,n}_k\left({\bf d}^m\right){\bf D}_0\left({\bf d}^m\right)+
\label{z26}\\
&&(\omega_q-k+1)\left(L^{q,n}_{k-1}\left({\bf d}^m\right){\bf G}_1\left({\bf 
d}^m\right)-M^{q,n}_{k-1}\left({\bf d}^m\right){\bf D}_1\left({\bf d}^m\right)
\right)-\nonumber\\
&&{\omega_q-k+2\choose 2}\left(L^{q,n}_{k-2}\left({\bf d}^m\right){\bf G}_2
\left({\bf d}^m\right)-M^{q,n}_{k-2}\left({\bf d}^m\right){\bf D}_2\left({\bf 
d}^m\right)\right)+\ldots-\nonumber\\
&&(-1)^k{\omega_q-2\choose k-2}\left(L^{q,n}_2\left({\bf d}^m\right){\bf G}_{
k-2}\left({\bf d}^m\right)-M^{q,n}_2\left({\bf d}^m\right){\bf D}_{k-2}\left(
{\bf d}^m\right)\right)+\nonumber\\
&&(-1)^k{\omega_q-1\choose k-1}\left(L^{q,n}_1\left({\bf d}^m\right){\bf G}_{
k-1}\left({\bf d}^m\right)-M^{q,n}_1\left({\bf d}^m\right){\bf D}_{k-1}\left(
{\bf d}^m\right)\right)=0\;,\hspace{1cm} 1\leq k\leq \omega_q\;,\nonumber
\eea
where $L^{q,n}_k\left({\bf d}^m\right)=M^{q,n}_k\left({\bf d}^m\right)=0$ if $k
\leq 0$.
\subsubsection{${\mathbb B}^{q,n}\left(s,{\bf d}^m\right)=0$}\label{s532}
In this section we find the 2nd system of linear equations (\ref{z33}) for 
$L_k^{q,n}\left({\bf d}^m\right)$ and $M_k^{q,n}\left({\bf d}^m\right)$. 
Making use of Reps (\ref{h1}) and (\ref{h2}) substitute Reps (\ref{g13}) and 
(\ref{g13a}) into expression for ${\mathbb B}^{q,n}\left(s,{\bf d}^m\right)$ 
given in section \ref{s53}. Thus, we present the 2nd Eq. in (\ref{z22a}) as 
follows,
\bea
&&M_1^{q,n}\left({\bf d}^m\right)\left[s^{\omega_q-1}-{\mathbb C}^{q,n}_{
\omega_q-1}\left(s,{\bf d}^m\right)\right]+M_2^{q,n}\left({\bf d}^m\right)
\left[s^{\omega_q-2}-{\mathbb C}^{q,n}_{\omega_q-2}\left(s,{\bf d}^m\right)
\right]+\ldots+\quad\label{z30a}\\
&&M_{\omega_q-1}^{q,n}\left({\bf d}^m\right)\left[s-{\mathbb C}^{q,n}_1\left(s,
{\bf d}^m\right)\right]+M_{\omega_q}^{q,n}\left({\bf d}^m\right)\left[1-
{\mathbb C}^{q,n}_0\left(s,{\bf d}^m\right)\right]-\nonumber\\
&&L_1^{q,n}\left({\bf d}^m\right){\mathbb S}^{q,n}_{\omega_q-1}\left(s,{\bf d}^m
\right)-L_2^{q,n}\left({\bf d}^m\right){\mathbb S}^{q,n}_{\omega_q-2}\left(
s,{\bf d}^m\right)-\ldots-L_{\omega_q}^{q,n}\left({\bf d}^m\right){\mathbb S}^{
q,n}_0\left(s,{\bf d}^m\right)=0\;.\hspace{1cm}\nonumber
\eea
Keeping in mind (\ref{z24}), after binomial expansion of the terms $\left(s-C_{
j,i}\right)^p$ in (\ref{h1}) and (\ref{h2}) and subsequent substitution into 
(\ref{z30a}) we recast the power terms in $s$ and get the final equality
\bea
&&\left[M^{q,n}_1\left({\bf d}^m\right)s^{\omega_q-1}+M^{q,n}_2\left({\bf d}^m
\right)s^{\omega_q-2}+\ldots+M^{q,n}_{\omega_q}\left({\bf d}^m\right)\right]
{\bf G}_0\left({\bf d}^m\right)+\label{z30b}\\
&&\left[M^{q,n}_1\left({\bf d}^m\right){\omega_q-1\choose 1}s^{\omega_q-2}+  
M^{q,n}_2\left({\bf d}^m\right){\omega_q-2\choose 1}s^{\omega_q-3}+\ldots+
M^{q,n}_{\omega_q-1}\left({\bf d}^m\right)\right]{\bf G}_1\left({\bf d}^m\right)
-\nonumber\\
&&\left[M^{q,n}_1\left({\bf d}^m\right){\omega_q-1\choose 2}s^{\omega_q-3}\!+
M^{q,n}_2\left({\bf d}^m\right){\omega_q-2\choose 2}s^{\omega_q-4}+\!\ldots\!+
M^{q,n}_{\omega_q-2}\left({\bf d}^m\right)\right]{\bf G}_2\left({\bf d}^m\right)
+\!\ldots\!-\nonumber\\
&&(-1)^{\omega_q}\left[M^{q,n}_1\left({\bf d}^m\right){\omega_q-1\choose
\omega_q-2}s+M^{q,n}_2\left({\bf d}^m\right)\right]{\bf G}_{\omega_q-2}\left(
{\bf d}^m\right)+
(-1)^{\omega_q}\;M^{q,n}_1\left({\bf d}^m\right){\bf G}_{\omega_q-1}\left(
{\bf d}^m\right)-\nonumber\\
&&\left[L^{q,n}_1\left({\bf d}^m\right)s^{\omega_q-1}+L^{q,n}_2\left({\bf d}^m
\right)s^{\omega_q-2}+\ldots+L^{q,n}_{\omega_q}\left({\bf d}^m\right)\right]
{\bf D}_0\left({\bf d}^m\right)+\nonumber\\
&&\left[L^{q,n}_1\left({\bf d}^m\right){\omega_q-1\choose 1}s^{\omega_q-2}+ 
L^{q,n}_2\left({\bf d}^m\right){\omega_q-2\choose 1}s^{\omega_q-3}+\ldots+
L^{q,n}_{\omega_q-1}\left({\bf d}^m\right)\right]{\bf D}_1\left({\bf d}^m\right)
-\nonumber\\
&&\left[L^{q,n}_1\left({\bf d}^m\right){\omega_q-1\choose 2}s^{\omega_q-3}\!+
L^{q,n}_2\left({\bf d}^m\right){\omega_q-2\choose 2}s^{\omega_q-4}+\!\ldots\!+
L^{q,n}_{\omega_q-2}\left({\bf d}^m\right)\right]{\bf D}_2\left({\bf d}^m\right)
+\!\ldots\!-\nonumber\\
&&(-1)^{\omega_q}\left[L^{q,n}_1\left({\bf d}^m\right)\!{\omega_q-1\choose
\omega_q-2}s+L^{q,n}_2\left({\bf d}^m\right)\right]{\bf D}_{\omega_q-2}\left(
{\bf d}^m\right)+
(-1)^{\omega_q}L^{q,n}_1\left({\bf d}^m\right){\bf D}_{\omega_q-1}\left(
{\bf d}^m\right)=0.\nonumber
\eea
Equating the corresponding contributions coming from the power terms $s^a$,   
$a=0,\ldots,\omega_q-1$, in the l.h.s. and the r.h.s. of (\ref{z30b}) we get
the 2nd system of $\omega_q$ linear equations for $L_k^{q,n}\left({\bf d}^m
\right)$ and $M_k^{q,n}\left({\bf d}^m\right)$. For short we present here the 
three first equations
\bea
&&a=\omega_q-1\hspace{2cm}
M^{q,n}_1\left({\bf d}^m\right){\bf G}_0\left({\bf d}^m\right)-
L^{q,n}_1\left({\bf d}^m\right){\bf D}_0 \left({\bf d}^m\right)=0\;,
\label{z30}\\\nonumber\\
&&a=\omega_q-2\hspace{2cm}
M^{q,n}_2\left({\bf d}^m\right){\bf G}_0\left({\bf d}^m\right)-
L^{q,n}_2\left({\bf d}^m\right){\bf D}_0\left({\bf d}^m\right)+\label{z30e}\\
&&\hspace{4cm}{\omega_q-1\choose 1}\left(M^{q,n}_1\left({\bf d}^m\right)
{\bf G}_1\left({\bf d}^m\right)+L^{q,n}_1\left({\bf d}^m\right)
{\bf D}_1\left({\bf d}^m\right)\right)=0\;,\nonumber\\\nonumber\\
&&a=\omega_q-3\hspace{2cm}
M^{q,n}_3\left({\bf d}^m\right){\bf G}_0\left({\bf d}^m\right)-
L^{q,n}_3\left({\bf d}^m\right){\bf D}_0\left({\bf d}^m\right)+\nonumber\\
&&\hspace{4cm}{\omega_q-2\choose 1}\left(M^{q,n}_2\left({\bf d}^m\right)
{\bf G}_1\left({\bf d}^m\right)+L^{q,n}_2\left({\bf d}^m\right)
{\bf D}_1\left({\bf d}^m\right)\right)-\nonumber\\
&&\hspace{4cm}{\omega_q-1\choose 2}\left(M^{q,n}_1\left({\bf d}^m\right)
{\bf G}_2\left({\bf d}^m\right)+L^{q,n}_1\left({\bf d}^m\right)
{\bf D}_2\left({\bf d}^m\right)\right)=0\;,\nonumber
\eea
and the last one
\bea
&&a=0\hspace{1.5cm}
M^{q,n}_{\omega_q}\left({\bf d}^m\right){\bf G}_0\left({\bf d}^m\right)
-L^{q,n}_{\omega_q}\left({\bf d}^m\right){\bf D}_0\left({\bf d}^m\right)+
\nonumber\\
&&\hspace{2.4cm}
M^{q,n}_{\omega_q-1}\left({\bf d}^m\right){\bf G}_1\left({\bf d}^m\right)+
L^{q,n}_{\omega_q-1}\left({\bf d}^m\right){\bf D}_1\left({\bf d}^m\right)-\ldots
+\nonumber\\
&&\hspace{2.5cm}(-1)^{\omega_q}\left(M^{q,n}_2\left({\bf d}^m\right)
{\bf G}_{\omega_q-2}\left({\bf d}^m\right)+L^{q,n}_2\left({\bf d}^m\right)
{\bf D}_{\omega_q-2}\left({\bf d}^m\right)\right)+\nonumber\\
&&\hspace{2.5cm}(-1)^{\omega_q}\left(M^{q,n}_1\left({\bf d}^m\right)
{\bf G}_{\omega_q-1}\left({\bf d}^m\right)+L^{q,n}_1\left({\bf d}^m\right)
{\bf D}_{\omega_q-1}\left({\bf d}^m\right)\right)=0\;.\nonumber
\eea
The whole set of these equations can be written in more compact form
\bea
&&M^{q,n}_k\left({\bf d}^m\right){\bf G}_0\left({\bf d}^m\right)-
L^{q,n}_k\left({\bf d}^m\right){\bf D}_0\left({\bf d}^m\right)+\label{z33}\\
&&(\omega_q-k+1)\left(M^{q,n}_{k-1}\left({\bf d}^m\right){\bf G}_1\left({\bf
d}^m\right)+L^{q,n}_{k-1}\left({\bf d}^m\right){\bf D}_1\left({\bf d}^m\right)
\right)-\nonumber\\
&&{\omega_q-k+2\choose 2}\left(M^{q,n}_{k-2}\left({\bf d}^m\right){\bf G}_2
\left({\bf d}^m\right)+L^{q,n}_{k-2}\left({\bf d}^m\right){\bf D}_2\left({\bf
d}^m\right)\right)+\ldots-\nonumber\\
&&(-1)^k{\omega_q-2\choose k-2}\left(M^{q,n}_2\left({\bf d}^m\right){\bf G}_{
k-2}\left({\bf d}^m\right)+L^{q,n}_2\left({\bf d}^m\right){\bf D}_{k-2}\left(
{\bf d}^m\right)\right)+\nonumber\\
&&(-1)^k{\omega_q-1\choose k-1}\left(M^{q,n}_1\left({\bf d}^m\right){\bf G}_{
k-1}\left({\bf d}^m\right)+L^{q,n}_1\left({\bf d}^m\right){\bf D}_{k-1}\left(
{\bf d}^m\right)\right)=0\;,\hspace{1cm} 1\leq k\leq \omega_q\;,\nonumber
\eea
where $L^{q,n}_k\left({\bf d}^m\right)=M^{q,n}_k\left({\bf d}^m\right)=0$ if $k
\leq 0$.
\subsubsection{Common Solutions of Equations ${\mathbb A}^{q,n}\left(s,{\bf d}^m
\right)=0$ and ${\mathbb B}^{q,n}\left(s,{\bf d}^m\right)=0$}\label{s533}
In this section we finish a proof of Theorem \ref{the2} stated in section 
\ref{s21}. For this purpose we combine two sets of linear equations (\ref{z26}) 
and (\ref{z33}) for ${\bf G}_r\left({\bf d}^m\right)$ and ${\bf D}_r\left({\bf 
d}^m\right)$ and solve them together. Start with two linear equations 
(\ref{z27}) and (\ref{z30}) for ${\bf G}_0\left({\bf d}^m\right)$ and ${\bf D}_0
\left({\bf d}^m\right)$,
\bea
L^{q,n}_1\left({\bf d}^m\right){\bf G}_0\left({\bf d}^m\right)+   
M^{q,n}_1\left({\bf d}^m\right){\bf D}_0 \left({\bf d}^m\right)=0,\quad\quad
M^{q,n}_1\left({\bf d}^m\right){\bf G}_0\left({\bf d}^m\right)-   
L^{q,n}_1\left({\bf d}^m\right){\bf D}_0 \left({\bf d}^m\right)=0,\nonumber
\eea
with discriminant $\Delta^{q,n}_1\left({\bf d}^m\right)=\left[L^{q,n}_1\left(
{\bf d}^m\right)\right]^2+\left[M^{q,n}_1\left({\bf d}^m\right)\right]^2$. 
Noting that by (\ref{y8}), or by Lemma \ref{lem1}, both $L^{q,n}_1\left({\bf d}
^m\right)$ and $M^{q,n}_1\left({\bf d}^m\right)$ do not vanish simultaneously, 
we arrive at trivial solution
\bea
{\bf G}_0\left({\bf d}^m\right)={\bf D}_0\left({\bf d}^m\right)=0\;.\label{z32}
\eea
Keeping in mind (\ref{z32}) continue with two equations (\ref{z27a}) and 
(\ref{z30e}) for ${\bf G}_1\left({\bf d}^m\right)$ and ${\bf D}_1\left({\bf d}^m
\right)$,
\bea
L^{q,n}_1\left({\bf d}^m\right){\bf G}_1\left({\bf d}^m\right)+
M^{q,n}_1\left({\bf d}^m\right){\bf D}_1 \left({\bf d}^m\right)=0,\quad\quad
M^{q,n}_1\left({\bf d}^m\right){\bf G}_1\left({\bf d}^m\right)-L^{q,n}_1\left(
{\bf d}^m\right){\bf D}_1 \left({\bf d}^m\right)=0,\nonumber
\eea
with the same discriminant $\Delta^{q,n}_1\left({\bf d}^m\right)\neq 0$ 
discussed above. By the same reason as for (\ref{z32}) we obtain ${\bf G}_1
\left({\bf d}^m\right)={\bf D}_1\left({\bf d}^m\right)=0$. Prove by induction 
that 
\bea
{\bf G}_r\left({\bf d}^m\right)={\bf D}_r\left({\bf d}^m\right)=0\;,\quad 1\leq 
r<\omega_q\;.\label{z34}
\eea
Indeed, let ${\bf G}_r\left({\bf d}^m\right)={\bf D}_r\left({\bf d}^m\right)=0$,
$1\leq r< r_*$ and write two linear equations (\ref{z26}) and (\ref{z33}) for 
$k=r_*+1$. Simplifying both equations we get
\bea
L^{q,n}_1\left({\bf d}^m\right){\bf G}_{r_*}\left({\bf d}^m\right)+
M^{q,n}_1\left({\bf d}^m\right){\bf D}_{r_*}\left({\bf d}^m\right)=0,\quad
M^{q,n}_1\left({\bf d}^m\right){\bf G}_{r_*}\left({\bf d}^m\right)-L^{q,n}
_1\left({\bf d}^m\right){\bf D}_{r_*}\left({\bf d}^m\right)=0,\nonumber
\eea
with discriminant $\Delta^{q,n}_1\left({\bf d}^m\right)\neq 0$. Thus, we arrive
at trivial solution ${\bf G}_{r_*}\left({\bf d}^m\right)={\bf D}_{r_*}\left(
{\bf d}^m\right)=0$, and therefore the existence of the general solutions 
(\ref{z34}) is proven.

On basis of ${\bf G}_r\left({\bf d}^m\right)$ and ${\bf D}_r\left({\bf d}^m
\right)$ build two complex expressions and make use of (\ref{z34}),
\bea
{\bf G}_0\left({\bf d}^m\right)-i\;{\bf D}_0\left({\bf d}^m\right)=0\;,\quad
\quad {\bf G}_r\left({\bf d}^m\right)+i\;{\bf D}_r\left({\bf d}^m\right)=0\;,
\quad 1\leq r<\omega_q\;.\label{z35}
\eea
By definition (\ref{z24}) write two equalities (\ref{z35}) in exponential form
and arrive at two equalities of Theorem \ref{the2}. Thus, this finishes proof 
of Theorem \ref{the2}.

We finish this section calculating the total number $N\left({\bf d}^m\right)$ of
polynomial and quasipolynomial identities imposed on the syzygies degrees of 
numerical semigroup ${\sf S}\left({\bf d}^m\right)$. By Theorem \ref{the1} we 
have $m$ polynomial identities including (\ref{bet2}). By Theorem \ref{the2} for
any pair $n$ and $q\geq 2$ we can write $\omega_q$ quasipolynomial identities 
such that $\omega_q=0$ if $q\nmid d_j$, $1\leq j\leq m$. Keeping in mind a range
of valuation for $n$, $\gcd(n,q)=1$ and $1\leq n<q/2$, we have for a fixed $q
\geq 3$ exactly $\phi(q)/2$ values of $n$, where $\phi(q)$ stands for the totient
function. The value $q=2$ is a special one and by (\ref{z18}) it corresponds to 
the value $n=1$ only. Since $q$ has to divide at least one generator $d_j$ of 
the tuple ${\bf d}^m=\left\{d_1,\ldots,d_m\right\}$ we arrive finally at
\bea
N\left({\bf d}^m\right)=m+\omega_2+\frac1{2}\sum_{q=3}^{\max{\bf 
d}^m}\omega_q\phi(q)\;.\label{z36}
\eea
\vspace{-1cm}
\section{Applications}\label{s6}
In this section we discuss different applications of Theorems \ref{the1} and 
\ref{the2} to the various kind of numerical semigroups and find more compact 
form of polynomial identities. We illustrate a validity of identities by 
examples for numerical semigroups discussed earlier in literature.
\subsection{Complete Intersections}\label{s61}
Start with a simple case of numerical semigroups where its Hilbert series is 
given by (\ref{teb4}) and a power sums of the degrees of syzygies read for $1
\leq k<m$,
\bea
\sum_{j=1}^{\beta_1}C_{j,1}^k=\sum_{j=1}^{m-1}e_j^k\;,\quad
\sum_{j=1}^{\beta_2}C_{j,2}^k=\sum_{j>k=1}^{m-1}(e_j+e_r)^k\;,\quad\ldots,\quad 
C_{1,m-1}^k=E^k\;,\quad E=\sum_{j=1}^{m-1}e_j\;.\label{z38}
\eea
Substituting (\ref{z38}) into (\ref{10c}) in Theorem \ref{the1} we get a set of 
algebraic identities for the tuple of integers $(e_1,\ldots,e_{m-1})$ when 
$1\leq k\leq m-2$,
\bea
\sum_{j=1}^{m-1}e_j^k-\sum_{j>k=1}^{m-1}(e_j+e_r)^k+\sum_{j>k>l=1}^{m-1}(e_j+
e_r+e_l)^k-\ldots+(-1)^m\left(\sum_{j=1}^{m-1}e_j\right)^k=0\;.\label{z40}
\eea
Identity (\ref{z40}) looks very {\em 'combinatorial'} and, indeed, for $k=1$ it 
can be verified trivially,
\bea
\sum_{j=1}^{m-1}e_j-\!\!\sum_{j>k=1}^{m-1}\!(e_j+e_r)+\ldots +(-1)^mE=\!
\left[1-\!{m-2\choose 1}\!+\ldots+(-1)^m\right]\!E=(1-1)^{m-2}E=0\;.
\nonumber
\eea
However, for $1<k<m-1$ such straightforward way of verification is hardly 
performable. This is why in Appendix \ref{appendix1} we give its verification 
based on the inclusion-exclusion principle applied to the set of tuples of 
length $k$ comprised elements of $m-1$ different sorts
\footnote{Combinatorial proof of (\ref{z40}) was kindly communicated to me by 
R. Pinchasi (Technion).}. 

What is much more important that this approach allows to calculate the left hand
side in (\ref{z40}) when $k=m-1$. Namely, the following identity holds,
\bea
\sum_{j=1}^{m-1}e_j^{m-1}-\sum_{j>k=1}^{m-1}(e_j+e_r)^{m-1}+\ldots+(-1)^m\left(
\sum_{j=1}^{m-1}e_j\right)^{m-1}\!\!=(-1)^m(m-1)!\prod_{i=1}^{m-1}e_i\;.
\label{z41}
\eea
Combining now the last identity with the 2nd identity (\ref{10d}) in Theorem 
\ref{the1} we get finally
\begin{corollary}\label{cor2}
Let the complete intersection semigroup ${\sf S}\left({\bf d}^m\right)$ be given
with its Hilbert series $H\left({\bf d}^m;z\right)$ in accordance with 
(\ref{tr16}) and (\ref{teb4}). Then the following identity holds,
\bea
\prod_{i=1}^{m-1}e_i=\prod_{i=1}^md_i\;.\label{w4}
\eea
\end{corollary}
In other words, the entire set of $m-1$ polynomial identities of Theorem 
\ref{the1} is reduced to only one nontrivial identity (\ref{w4}). Regarding 
Theorem \ref{the2}, its quasipolynomial relations imply another set of 
identities which definitely do not have combinatorial origin.
\begin{corollary}\label{cor3}
Let the complete intersection semigroup ${\sf S}\left({\bf d}^m\right)$ be given
with its Hilbert series $H\left({\bf d}^m;z\right)$ in accordance with 
(\ref{tr16}) and (\ref{teb4}). Then for every $1<q\leq\max\left\{d_1,\ldots,d_m
\right\}$, and $\gcd(n,q)=1$, $1\leq n<q/2$, the following quasipolynomial 
identities hold,
{\footnotesize
\bea
\sum_{j=1}^{m-1}\!e_j^r\exp\left(i\frac{2\pi n}{q}e_j\right)-\!\sum_{j>l=1}^{m-
1}\!\!\left(e_j+e_l\right)^r\exp\left(i\frac{2\pi n}{q}\left(e_j+e_l\right)
\right)+\ldots+(-1)^m\!\sum_{j=1}^{m-1}\!\left(E-e_j\right)^r\exp\left(i
\frac{2\pi n}{q}\left(E-e_j\right)\right)\!=0.\nonumber
\eea}
where $1\leq r<\omega_q$ and $E=\sum_{i=1}^{m-1}\!e_i$. However, in the case 
$r=0$ another identity holds,
{\footnotesize
\bea
\sum_{j=1}^{m-1}\!\exp\left(i\frac{2\pi n}{q}e_j\right)-\!\sum_{j>l=1}^{m-1}\!
\exp\left(i\frac{2\pi n}{q}\left(e_j+e_l\right)\right)+\ldots+(-1)^m\!
\sum_{j=1}^{m-1}\!\exp\left(i\frac{2\pi n}{q}\left(E-e_j\right)\right)\!=1\;.
\nonumber
\eea}
\end{corollary}
Below we illustrate Corollaries \ref{cor2} and \ref{cor3} by presenting the 
Hilbert series of two complete intersection semigroups generated by four 
elements \cite{de76}, \cite{kra85}: a) $(8,9,10,12)$ and b) $(10,14,15,21)$,
\bea
H_a(z)=\frac{(1-z^{18})(1-z^{20})(1-z^{24})}{\left(1-z^8\right)\left(1-z^9
\right)\left(1-z^{10}\right)\left(1-z^{12}\right)},\quad
H_b(z)=\frac{(1-z^{30})(1-z^{35})(1-z^{42})}{\left(1-z^{10}\right)\left(1-z^{14}
\right)\left(1-z^{15}\right)\left(1-z^{21}\right)}\;.\nonumber
\eea
For both semigroups the identities of both Corollaries \ref{cor2} and \ref{cor3}
are satisfied.
\subsubsection{Telescopic Semigroups}\label{s611}
Two semigroups ${\sf S}\left(\{8,9,10,12\}\right)$ and ${\sf S}\left(\{10,14,15,
21\}\right)$ discussed in the previous section \ref{s61} present two different 
kinds of complete intersections: telescopic and complete intersections-not 
telescopic semigroups, respectively. The semigroups of the 1st kind present the 
most simple sort of complete intersections. 

Following \cite{kp95} start with definition. For given numerical semigroup 
${\sf S}\left({\bf d}^m\right)$ with generators ${\bf d}^m=\{d_1,\ldots,d_m\}$, 
$\gcd(d_1,\ldots,d_m)=1$, (not necessarily in increasing order), let us denote 
$g_k=\gcd(d_1,\ldots,d_k)>1$ and ${\sf S}_k={\sf S}\left(\left\{\frac{d_1}{g_k},
\frac{d_2}{g_k},\ldots,\frac{d_k}{g_k}\right\}\right)$ for $1\leq k<m$, and $g_1
=d_1$, $g_m=1$. Then ${\sf S}\left({\bf d}^m\right)$ is said to be telescopic 
iff $\frac{d_k}{g_k}\in{\sf S}_{k-1}$ for all $2\leq k\leq m$. Regarding its 
Hilbert series (\ref{teb4}), the syzygy degrees $e_j$ of the 1st kind follow by 
calculating $m-1$ minimal linear relations between the generators $d_j$. Their 
straightforward calculation gives,
\bea
e_1=lcm(d_1,d_2)\;,\quad e_2=d_3\frac{g_2}{g_3}\;,\ldots,\quad e_j=d_{j+1}
\frac{g_j}{g_{j+1}}\;,\quad 2\leq j<m\;,\quad e_{m-1}=d_m\;g_{m-1}\;.\quad
\label{w5w}
\eea
As for the polynomial identity (\ref{w4}), it becomes trivial in the case of
telescopic semigroups. Indeed, this can be verified if substituting (\ref{w5w}) 
into (\ref{w4}). Thus, there remain only the set of quasipolynomial identities 
after substitution (\ref{w4}) into identities of Corollary \ref{cor3}.

\subsection{Symmetric Numerical Semigroups}\label{s62}
This kind of semigroups is of high importance due to the theorem of Kunz 
\cite{ku70} which asserts that a graded polynomial subring, associated with 
semigroup ${\sf S}\left({\bf d}^m\right)$, is Gorenstein iff ${\sf S}\left({\bf 
d}^m\right)$ is symmetric. For short, denote ${\cal Q}\left({\bf d}^m\right)=
\deg Q\left({\bf d}^m;z\right)$ and define the following combination of Betti 
numbers,
\bea
{\cal B}\left({\bf d}^m\right)=-1+\beta_1\left({\bf d}^m\right)-\beta_2\left(
{\bf d}^m\right)+\ldots +(-1)^{\mu}\beta_{\mu-1}\left({\bf d}^m\right)\;,\quad
\mu=\left\lfloor \frac{m}{2}\right\rfloor\;.\label{w5c}
\eea
Then, by Theorem \ref{the1} for symmetric semigroups with even edim $2m$ it 
holds for $1\leq k\leq 2m-1$,
\bea
{\cal Q}^k\left({\bf d}^{2m}\right)-\sum_{r=1}^{m-1}(-1)^r\sum_{j=1}^{\beta_r
\left({\bf d}^{2m}\right)}\left\{C_{j,r}^k-\left[{\cal Q}\left({\bf d}^{2m}
\right)-C_{j,r}\right]^k\right\}=(2m-1)!\;\pi_{2m}\;\delta_{k,2m-1}\;.\label{w5}
\eea
In the case $k=1$ the last identity reads
gives
\bea
\sum_{j=1}^{\beta_1\left({\bf d}^{2m}\right)}\!\!\!\!C_{j,1}-\sum_{j=1}^{
\beta_2\left({\bf d}^{2m}\right)}\!\!\!\!C_{j,2}+\ldots+(-1)^m\sum_{j=1}^{
\beta_{m-1}\left({\bf d}^{2m}\right)}\!\!\!\!C_{j,m-1}=\frac1{2}{\cal B}\left(
{\bf d}^{2m}\right){\cal Q}\left({\bf d}^{2m}\right)\;.\label{w5a}
\eea
By the same Theorem \ref{the1} for symmetric semigroups with odd edim=$2m+1$ it 
holds for $1\leq k\leq 2m$,
\bea
{\cal Q}^k\left({\bf d}^{2m+1}\right)+\sum_{r=1}^{m-1}(-1)^r\!\!\!\!
\sum_{j=1}^{\beta_r\left({\bf d}^{2m+1}\right)}\!\!\!\!
\left\{C_{j,r}^k+\left[{\cal Q}\left({\bf d}^{2m+1}\right)-C_{j,r}\right]^k
\right\}+(-1)^m\!\!\!\!\sum_{j=1}^{\beta_m\left({\bf d}^{2m+1}\right)}\!\!\!\!
C_{j,m}^k+\nonumber\\(2m)!\;\pi_{2m+1}\;\delta_{k,2m}=0\;,\label{w5d}
\eea
that in the case $k=1$ gives
\bea
\sum_{j=1}^{\beta_m\left({\bf d}^{2m+1}\right)}\!\!\!\!C_{j,m}={\cal B}\left(
{\bf d}^{2m+1}\right){\cal Q}\left({\bf d}^{2m+1}\right).\label{w5b}
\eea
\subsubsection{Symmetric Semigroups ${\sf S}\left({\bf d}^4\right)$ (not 
complete intersections)}\label{s621}
In 1975, Bresinsky \cite{ber752} has shown that symmetric semigroups ${\sf S}
\left({\bf d}^4\right)$, which are not complete intersections, have always 
$\beta_1\left({\bf d}^4\right)=5$. Denoting invariants $I_k=\sum_{j=1}^5C_{j,1}
^k$, we obtain by (\ref{w5}) two different polynomial identities,
\bea
8I_3-6I_2I_1+I_1^3=24\pi_4\;,\quad {\cal Q}\left({\bf d}^4\right)=\frac1{2}\;
I_1\;.\label{w6}
\eea
Below we illustrate identities (\ref{w6}) by presenting the Hilbert series of 
two symmetric semigroups (not complete intersection) generated by four elements
\cite{de76},
\bea
H\left(\{5,6,7,8\};z\right)=\frac{1-z^{12}-z^{13}-z^{14}-z^{15}-z^{16}+z^{19}
+z^{20}+z^{21}+z^{22}+z^{23}-z^{35}}{\left(1-z^5\right)\left(1-z^6\right)
\left(1-z^7\right)\left(1-z^8\right)}\;,\nonumber
\eea
\bea
H\left(\{8,13,15,17\};z\right)=\frac{1-z^{30}-z^{32}-z^{34}-z^{39}-z^{41}+
z^{47}+z^{49}+z^{54}+z^{56}+z^{58}-z^{88}}{\left(1-z^8\right)\left(1-z^{13}
\right)\left(1-z^{15}\right)\left(1-z^{17}\right)}\;.\nonumber
\eea
For both semigroups the quasipolynomial identities of Theorem \ref{the2} are 
also satisfied.
\subsubsection{Symmetric Semigroups ${\sf S}\left({\bf d}^5\right)$ (not 
complete intersections)}\label{s622}
In contrast to symmetric semigroups ${\sf S}\left({\bf d}^4\right)$, the problem
of admissible values of the 1st Betti number for symmetric semigroups (not 
complete intersections) ${\sf S}\left({\bf d}^m\right)$, $m\geq 5$, or their 
upper bounds, are still open (see \cite{bar06}, Problem 3). Therefore in this 
section we study the polynomial identities of Theorem \ref{the1} for $m=5$ and 
arbitrary $\beta_1\left({\bf d}^5\right)$. 

Henceforth, we skip for short the notation ${\bf d}^5$ in the Betti numbers and 
denote symmetric invariants $J_{r,k}=\sum_{j=1}^{\beta_r}C_{j,r}^k$ where $r=1,
2$ and $1\leq k\leq 4$. By (\ref{w5d}) and (\ref{w5b}) we get four different 
polynomial identities,
\bea
J_{2,1}=(\beta_1-1){\cal Q}\left({\bf d}^5\right),\hspace{1cm}
J_{2,1}(2J_{1,1}-J_{2,1})+(\beta_1-1)\left(J_{2,2}-2J_{1,2}\right)=0,
\hspace{1.8cm}\label{z43}\\
J_{2,1}^2(3J_{1,1}-J_{2,1})-3(\beta_1-1)J_{1,2}J_{2,1}+(\beta_1-1)^2J_{2,3}=0,
\hspace{1cm}\nonumber\\
J_{2,1}^3(4J_{1,1}-J_{2,1})-6(\beta_1-1)J_{1,2}J_{2,1}^2+
4(\beta_1-1)^2J_{1,3}J_{2,1}+(\beta_1-1)^3\left(J_{2,4}-2J_{1,4}-24\pi_5\right)=
0.\nonumber
\eea
For verification we have chosen the known symmetric semigroup ${\sf S}\left(
\{19,23,29,31,37\}\right)$ generated by five elements and found by computer 
calculations by Bresinsky \cite{ber79}. Its corresponding Betti numbers read: 
$\beta_1=13,\;\beta_2=24,\;\beta_3=13,\;\beta_4=1$. The following numerator of 
its Hilbert series was calculated by the author,
\bea
Q(z)&=&1-z^{60}-z^{69}-z^{75}-z^{77}-z^{81}-z^{85}-z^{87}-z^{93}-z^{95}+z^{98}
-z^{99}+z^{100}-z^{103}+\nonumber\\
&&z^{104}-z^{105}+2z^{106}+z^{108}+z^{110}-z^{111}+z^{112}+z^{114}+z^{116}+
2z^{118}+2z^{122}+\nonumber\\
&&z^{124}+z^{126}+z^{128}-z^{129}+z^{130}+z^{132}+2z^{134}-z^{135}+z^{136}-
z^{137}+z^{140}-z^{141}+\nonumber\\
&&z^{142}-z^{145}-z^{147}-z^{153}-z^{155}-z^{159}-z^{163}-z^{165}-z^{171}-
z^{180}+z^{240}\;.\nonumber
\eea
Straightforward calculation shows that the polynomial identities (\ref{z43})
are satisfied as well as the quasipolynomial identities of Theorem \ref{the2}.
\subsection{Nonsymmetric Numerical Semigroups ${\sf S}\left({\bf d}^3\right)$}
\label{s63}
In 2004, studying the Frobenius problem for numerical semigroups ${\sf S}\left(
{\bf d}^3\right)$ with the Hilbert series
\bea
H\left({\bf d}^3;z\right)=\frac{1-z^{e_1}-z^{e_2}-z^{e_3}+z^{q_1}+z^{q_2}}
{\left(1-z^{d_1}\right)\left(1-z^{d_2}\right)\left(1-z^{d_3}\right)}\;,
\label{w8}
\eea
in \cite{ros04} and \cite{fel04} there were found the degrees $q_1$ and $q_2$, 
$q_1<q_2$, of syzygies of the 2nd kind in terms of 3 degrees $e_1$, $e_2$ and 
$e_3$ of syzygies of the 1st kind and 3 generators $d_1,d_2,d_3$,
\bea
q_{1,2}=\frac1{2}\left[(e_1+e_2+e_3)\pm\sqrt{e_1^2+e_2^2+e_3^2-2(e_1e_2+e_2e_3+
e_3e_1)+4d_1d_2d_3}\;\right]\;.\label{w9}
\eea
The last formula gives  by (\ref{bet1a}) the Frobenius number $F\left({\bf d}^3
\right)=q_2-\sigma_1$. Recently a derivation of formula (\ref{w9}) has been
shorten essentially \cite{aic09} by making use of the Ap\'ery set of ${\sf S}
\left({\bf d}^3\right)$ and of relation between its generating function and the 
Hilbert series. It turns out that the way to derive (\ref{w9}) can be shorten 
much more. Indeed, by Theorem \ref{the1} two identities hold
\bea
e_1+e_2+e_3=q_1+q_2\;,\quad\quad e_1^2+e_2^2+e_3^2=q_1^2+q_2^2-2d_1d_2d_3\;,
\label{w10}
\eea
and after trivial algebra we arrive at (\ref{w9}).

The derivation of the degrees $e_1$, $e_2$ and $e_3$ through the three 
generators $d_i$ is much more difficult problem which encounters the Curtis 
theorem \cite{curt90} on non algebraic Rep of the Frobenius number $F\left({\bf 
d}^3\right)$.  By (\ref{w9}) this statement is equivalent to the claim that 
neither of the syzygy degrees $e_i$ is representable by algebraic function in 
$d_1,d_2,d_3$. The analytic Rep for $e_i$ by integration in the complex plane 
was found in \cite{fel08}. However, based on quasipolynomial identities of 
Theorem \ref{the2} one can build another set of non algebraic Reps. E.g., in the
case when all generators $d_i$ are primes the syzygy degrees $e_i$ come as 
solutions of polynomial equations (\ref{w10}) together with exponential 
equations
\bea
&&\xi_{d_j}^{ne_1}+\xi_{d_j}^{ne_2}+\xi_{d_j}^{ne_3}=\xi_{d_j}^{nq_1}+
\xi_{d_j}^{nq_2}+1\;,\quad j=1,2,3\;,\quad \gcd(n,d_j)=1\;.\label{w11}
\eea
\subsubsection{Pseudosymmetric Semigroups ${\sf S}\left({\bf d}^3\right)$}
\label{s631}
A numerical semigroup ${\sf S}\left({\bf d}^m\right)$ is pseudosymmetric if 
$F\left({\bf d}^m\right)$ is even and the only integer such $s\in {\mathbb N}
\setminus{\sf S}\left({\bf d}^m\right)$ and $F\left({\bf d}^m\right)-s\not\in 
{\sf S}\left({\bf d}^m\right)$ is $s=1/2\;F\left({\bf d}^m\right)$. A case $m=3$
is most simple and allows to calculate the Frobenius number \cite{ros05}, 
\cite{f10},
\bea
F\left({\bf d}^3\right)=-\sigma_1+\sqrt{\sigma_1^2-4(d_1d_2+d_2d_3+d_3d_1)+
4d_1d_2d_3}\;,\label{x1}
\eea
and the whole numerator of the Hilbert series \cite{f10},
\bea
Q\left({\bf d}^3;z\right)=1-z^{d_1+d_2+\frac1{2}F\left({\bf d}^3\right)}-
z^{d_2+d_3+\frac1{2}F\left({\bf d}^3\right)}-z^{d_3+d_1+\frac1{2}F\left(
{\bf d}^3\right)}+z^{\frac1{2}F\left({\bf d}^3\right)+\sigma_1}+
z^{F\left({\bf d}^3\right)+\sigma_1}\;.\nonumber
\eea
The degrees of syzygies of the above expression satisfy the 1st Eq. in 
(\ref{w10}). As for the 2nd Eq. in (\ref{w10}), it is reduced to quadratic Eq. 
in $F\left({\bf d}^3\right)$ and gives (\ref{x1}).
\subsection{Numerical Semigroups of Maximal edim}\label{s64}
A semigroup of maximal edim (for short, MED semigroup), which is generated by 
tuple ${\bf d}^m_{MED}=\{m,d_2,\ldots,d_m\}$, is never symmetric. Many explicit 
results are known about its Betti numbers, genus and Frobenius number 
\cite{heku71}, \cite{sal79}, \cite{fa09}, e.g., $F\left({\bf d}^m_{MED}\right)=
d_m-m$. Regarding the Hilbert series, which was found recently \cite{f10}, the 
partial contributions $Q_k\left({\bf d}^m_{MED};z\right)$ (see (\ref{bet05})) 
to the whole numerator $Q\left({\bf d}^m_{MED};z\right)$ in Hilbert series 
(\ref{tr16}) read,
\bea
Q_k\left({\bf d}^m_{MED};z\right)=I_{m,k}(z)+J_{m,k}(z)\;,\quad 1\leq k\leq m-2
\;,\quad Q_{m-1}\left({\bf d}^m_{MED};z\right)=I_{m,m-1}(z)\;,\quad\label{w20}
\eea
where for $d_j\in\{d_2,\ldots,d_m\}$ the following notations stand,
\bea
I_{m,k}(z)=\!\!\!\!\sum_{j_1\ne j_2>\ldots >j_{k-1}\geq 2}^m\!\!\!\!\!\!\!
z^{2d_{j_1}+\overbrace{d_{j_2}+\ldots+d_{j_{k}}}^{{\rm k-1}\;terms}}\;,\quad 
J_{m,k}(z)=k\!\!\!\!\sum_{j_1>\ldots >j_k\geq 2}^m\!\!\!\!\!\!z^{\overbrace{
d_{j_1}+d_{j_2}+\ldots+d_{j_{k+1}}}^{{\rm k+1}\;terms}}.\quad\label{w21}
\eea
We give the polynomials $I_{m,k}(z)$ and $J_{m,k}(z)$ for small and large
indices $k$,
\bea
&&I_{m,1}(z)=\sum_{j_1\geq 2}^mz^{2d_{j_1}},\;\;I_{m,2}(z)=\sum_{j_1\neq j_2
\geq 2}^mz^{2d_{j_1}+d_{j_2}},\;\;I_{m,3}(z)=\sum_{j_1\neq j_2>j_3\geq 2}^m
z^{2d_{j_1}+d_{j_2}+d_{j_3}},\;\ldots ,\nonumber\\
&&I_{m,m-2}(z)=\sum_{j_1\neq j_2>\ldots >j_{m-2}\geq 2}^mz^{2d_{j_1}+d_{j_2}+
\ldots+ d_{j_{m-2}}},\;\;\;\;I_{m,m-1}(z)=z^{\sigma_1}\sum_{j\geq 2}^m
z^{d_j-m}\;,\label{w22}
\eea
\bea
J_{m,1}(z)=\!\!\sum_{j_1>j_2\geq 2}^m\!\!z^{d_{j_1}+d_{j_2}},\;J_{m,2}(z)=2
\!\!\sum_{j_1>j_2>j_3\geq 2}^m\!\!z^{d_{j_1}+d_{j_2}+d_{j_3}},\ldots ,\;
J_{m,m-2}(z)=(m-2)z^{\sigma_1-m}.\nonumber
\eea
Substituting the degrees of monomials (\ref{w20}) into two identities of 
Theorem \ref{the1} in accordance with an ordinal number $k$ of syzygies, which 
monomials (\ref{w20}) belong to, we obtain Corollary.
\begin{corollary}\label{cor4}
Let a set $\{d_2,\ldots,d_m\}$ of $m-1$ distinct positive integers $d_j$ be 
given. Then for $1\leq k\leq m-2$ the following identities hold,
\bea
2^k\sum_{j_1\geq 2}^md_{j_1}^k+\!\!\sum_{j_1>j_2\geq 2}^m\!\!\left(d_{j_1}+
d_{j_2}\right)^k-\!\!\sum_{j_1\neq j_2\geq 2}^m\!\!\left(2d_{j_1}+d_{j_2}
\right)^k-2\!\!\!\sum_{j_1>j_2>j_3\geq 2}^m\!\!\!\left(d_{j_1}+d_{j_2}+
d_{j_3}\right)^k+\quad\nonumber\\
\ldots+(-1)^m\sum_{j\geq 2}\left(\sigma_1+d_j-m\right)^k=0\;,\nonumber
\eea
where $\sigma_1=m+\sum_{j=2}^md_j$. If $k=m-1$, then
\bea
2^{m-1}\sum_{j_1\geq 2}^md_{j_1}^{m-1}+\!\!\sum_{j_1>j_2\geq 2}^m\!\!\left(
d_{j_1}+d_{j_2}\right)^{m-1}-\!\!\sum_{j_1\neq j_2\geq 2}^m\!\!\left(
2d_{j_1}+d_{j_2}\right)^{m-1}-\hspace{2cm}\nonumber\\
2\!\!\!\sum_{j_1>j_2>j_3\geq 2}^m\!\!\!
\left(d_{j_1}+d_{j_2}+d_{j_3}\right)^{m-1}
+\ldots+(-1)^m\sum_{j\geq 2}\left(\sigma_1+d_j-m\right)^{m-1}=(-1)^mm!
\prod_{j=2}^md_j\;.\nonumber
\eea
\end{corollary}
As an example of the MED semigroup one can take the semigroup ${\sf S}\left(
\{3,5,7\}\right)$ with Hilbert series given in (\ref{ar10b}). A straightforward 
calculation shows that the polynomial identities of Corollary \ref{cor4} are 
satisfied as well as the quasipolynomial identities of Theorem \ref{the2}.
\subsubsection{Almost Symmetric Semigroups of Maximal edim}\label{s641}
Almost symmetric semigroups ${\sf S}\left({\bf d}^m\right)$ were introduced in 
\cite{bafr97} as a generalization of the symmetric and pseudosymmetric ones. 
They can be defined \cite{bafr97} by equality $\beta_{m-1}\left({\bf d}^m\right)
=1+\#\Delta_{{\cal H}}\left({\bf d}^m\right)$, where $\Delta_{{\cal H}}\left(
{\bf d}^m\right)=\left\{h\not\in {\sf S}\left({\bf d}^m\right)\;|\;F\left({\bf 
d}^m\right)-h\not\in {\sf S}\left({\bf d}^m\right)\right\}$ is a subset of the 
set of the gaps $\Delta\left({\bf d}^m\right)$. A detailed study of such 
semigroups can be found in \cite{bafr97} and \cite{f10}. Here we discuss a 
special class of almost symmetric MED semigroups. Necessary and sufficient 
conditions for a minimal set ${\bf d}^m$ to generate such semigroup is that for 
every element $d_j\in {\bf d}^m$ there exists its counterpartner $d_{m-j+1}$ 
such that (see \cite{f10}, Theorem 7),
\bea
d_j+d_{m-j+1}=\frac{2\sigma_1}{m}\;,\;\;\;\;\;1\leq j\leq m\;.\label{h5}
\eea
Relations (\ref{h5}) do not reduce the total number of the polynomial identities
in Corollary \ref{cor4}. However, they do simplify a large number of terms 
contributing to these identities.

Below we illustrate the polynomial identities of Corollary \ref{cor4} by 
presenting the Hilbert series of almost symmetric MED semigroups ${\sf S}\left(
\{4,10,19,25\}\right)$ taken from \cite{bafr97}. Its corresponding Betti numbers
read: $\beta_1=6$, $\beta_2=8$, $\beta_3=3$, while the numerator of its Hilbert 
series was calculated in \cite{f10},
{\footnotesize
\bea
H(z)=\frac{1-z^{20}-z^{29}-z^{35}-z^{38}+z^{39}-z^{44}+z^{45}+z^{48}-z^{50}
+2z^{54}+z^{60}+z^{63}-z^{64}+z^{69}-z^{73}-z^{79}}{\left(1-z^4\right)
\left(1-z^{10}\right)\left(1-z^{19}\right)\left(1-z^{25}\right)}\;.\nonumber
\eea}
Straightforward calculation shows that the polynomial identities in Corollary 
\ref{cor4} are satisfied as well as the quasipolynomial identities of Theorem 
\ref{the2}.
\appendix
\renewcommand{\theequation}{\thesection\arabic{equation}}
\section{Appendix: Combinatorial proof of identities (\ref{z40}) and 
(\ref{z41})}\label{appendix1}
\setcounter{equation}{0}
Consider $n$ sets of elements $a_{j,r}\in {\mathbb E}_r$, $\# {\mathbb E}_r=e_
r$, $1\leq r\leq n$, of different colors and construct a tuple ${\bf a}^k$ of 
length $k$ which is composed of elements $a_{j,r}$. Denote by ${\mathbb C}_i$ 
a set of all ${\bf a}^k$-tuples which do not contain any element of the $i$th 
color, i.e., 
\bea
{\mathbb C}_i=\{{\bf a}^k\;|\;a_{j,r}\in {\bf a}^k,\;r\neq i\}\;,\quad
\#{\mathbb C}_i=\left(E-e_i\right)^k\;,\quad E=\sum_{j=1}^ne_j\;.\label{ap2}
\eea
Consider an intersection of two such sets ${\mathbb C}_i$ and ${\mathbb C}_l$ 
which do not contain any element of the $i$th and $l$th colors,
\bea
{\mathbb C}_i\cap {\mathbb C}_l=\{{\bf a}^k\;|\;a_{j,r}\in {\bf a}^k,\;r\neq i,
\;r\neq l\}\;,\quad \#\left\{{\mathbb C}_i\cap {\mathbb C}_l\right\}=
\left(E-e_i-e_l\right)^k\;.\label{ap3}
\eea
Continuing to build intersections of three and more sets ${\mathbb C}_i$ write 
by the inclusion-exclusion principle an identity for their cardinalities,
\bea
\#\bigcup_{i=1}^n{\mathbb C}_i=\sum_{i=1}^n\left(E-e_i\right)^k-\!\sum_{i>l=1}^n
\!\!\left(E-e_i-e_l\right)^k+\ldots-(-1)^n\sum_{i>l=1}^n\!\!\left(e_i+e_l\right)
^k+(-1)^n\sum_{i=1}^ne_i^k,\;\label{ap4}
\eea
where $\bigcup_{i=1}^n{\mathbb C}_i=\{{\bf a}^k\;|\;a_{j,r}\in{\bf a}^k,\;1
\leq r\leq n\}$. If $k<n$ then $\#\bigcup_{i=1}^n{\mathbb C}_i=E^k$ and 
therefore
\bea
\sum_{i=1}^ne_i^k-\sum_{i>l=1}^n\!\!\left(e_i+e_l\right)^k+\ldots+(-1)^n
\sum_{i=1}^n\left(E-e_i\right)^k -(-1)^nE^k=0\;.\label{ap5}
\eea
In the case $k=n$ we have $\#\bigcup_{i=1}^n{\mathbb C}_i=E^n-n!\prod_{i=1}^n
e_i$. Inserting the last into (\ref{ap4}) we get 
\bea
\sum_{i=1}^ne_i^n-\sum_{i>l=1}^n\!\!\left(e_i+e_l\right)^n+\ldots+(-1)^n\sum_{i
=1}^n\left(E-e_i\right)^n-(-1)^nE^n=-(-1)^nn!\prod_{i=1}^ne_i\;.\label{ap6}
\eea
Substituting $n=m-1$ into (\ref{ap5}) and (\ref{ap6}) we arrive at (\ref{z40}) 
and (\ref{z41}), respectively.
\section*{Acknowledgement}
The author thanks A. Juhasz, R. Pinchasi and B. Rubinstein for useful comments. 
The research was partly supported by the Kamea Fellowship.

\end{document}